\documentclass[12pt]{article}
\usepackage{amsfonts}
\usepackage{amsmath}
\usepackage{amssymb}
\usepackage{graphicx}

\usepackage{amsfonts}
\usepackage{amsmath}
\newtheorem{theorem}{Theorem}[section]
\newtheorem{lemma}[theorem]{Lemma}
\newtheorem{corollary}[theorem]{Corollary}

\newtheorem{example}[theorem]{Example}
\newtheorem{examples}[theorem]{Examples}

\newcommand{\dis}{\displaystyle}
\newcommand{\HS}[3]{\left(\frac{#1,#2}{#3}\right)}
\newcommand{\ds}{\displaystyle}

\newcommand{\IR}{\mathbb R}
\newcommand{\IQ}{\mathbb Q}
\newcommand{\IZ}{\mathbb Z}
\newcommand{\IS}{\mathbb S}
\newcommand{\IC}{\mathbb C}
\newcommand{\IH}{\mathbb H}

\newcommand{\G}{\Gamma}

\def\tr{\mbox{\rm{tr\,}}}

\def\PSL{\mbox{\rm{PSL}}}
\def\SL{\mbox{\rm{SL}}}

\def\log{\mbox{\rm{log}}}
\def\Isom{\mbox{\rm{Isom}}}

\def\Isom{\mbox{\rm{Isom}}}

\def\Ram{\mbox{\rm{Ram}}}

\def\em{\it}

\title{The $(p,q)$-arithmetic hyperbolic lattices; \\ $p,q \geq 6$}
\author{C. Maclachlan and G.J. Martin  \thanks{Research supported in
part by grants from the N.Z.
Marsden Fund and the New Zealand Royal Society (James Cook Fellowship). \newline \newline AMS
(1991) Classification.
Primary 30F40, 30D50, 20H10, 22E40, 53A35, 57M60}  }
\date{}
\begin{document}

\maketitle

\begin{abstract}
We prove there are exactly $16$ arithmetic lattices of hyperbolic $3$-space which are generated by two elements of finite orders $p$ and $q$ with $p,q\geq6$.    We also verify a conjecture of H.M. Hilden, M.T. Lozano, and J.M. 
Montesinos  concerning the orders of the singular sets of arithmetic orbifold Dehn surgeries on two bridge knot and link complements.
\end{abstract}

\section{Introduction}

There are infinitely many lattices in the group $\Isom^+\IH^3$ $\cong \PSL(2,
\IC)$, of orientation-preserving isometries of hyperbolic 3-space (equivalently
Kleinian groups of finite co-volume) which can be generated by two elements
of finite orders $p$ and $q$. For instance, all but finitely many $(p,0)-(q,0)$
Dehn surgeries on any of the infinitely many hyperbolic two-bridge links 
 will have fundamental groups which are such 
uniform (co-compact) lattices \cite{Thurston}. Two infinite families of such groups are shown below  in Figure 1. 

In \cite{MM}, we showed that, up to conjugacy,
only finitely many of these lattices can be arithmetic. In \cite{MM2}, we identified the 20 
such non-uniform lattices of which 15 were {\em generalised triangle groups};
that is, groups with a presentation of the form $\langle x,y : x^p = y^q = 
w(x,y)^r = 1 \rangle$ where $w(x,y)$ is a word involving both $x$ and $y$ (see 
\cite{FR1,BMS}) and $p,q,r\geq 2$.

\scalebox{0.5}{\includegraphics*[viewport=-50 500 700 800]{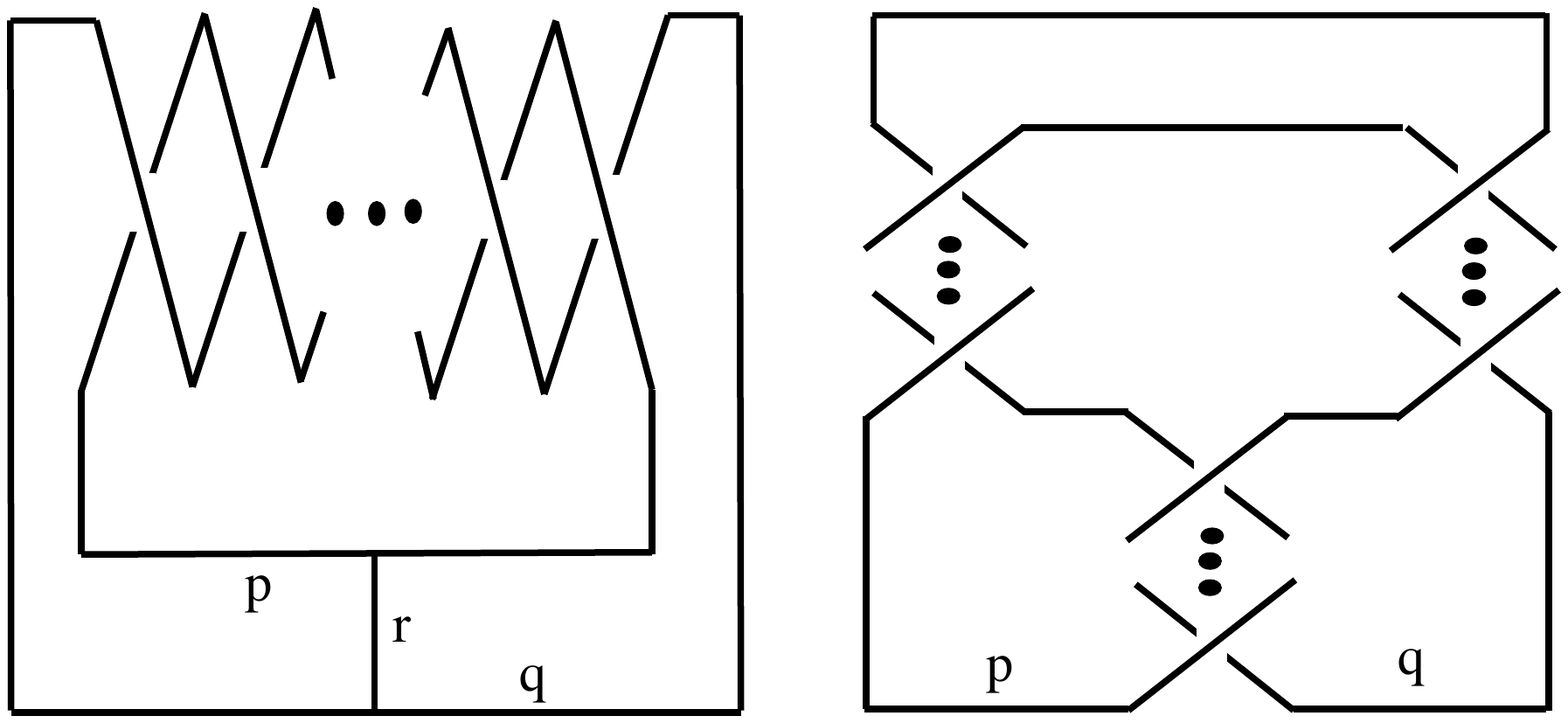}}
{\small Figure 1.   On the left the singular set in $\IS^3$ of an infinite family of hyperbolic generalised triangle groups (eg. $r=2$ and $p,q\geq 3$) with $p=q$ if there are an odd number of twists,  \cite{JR}.  On the right an infinite family of  hyperbolic two bridge $p,q$ links \cite{Thurston2}.}

 \bigskip

In this paper we prove that, up to conjugacy, there are exactly $16$
arithmetic lattices in $\Isom^+\IH^3$ which can be generated by two elements 
of finite orders $p$ and $q$ with $6 \leq p,q $. Among these groups there are,
curiously, no generalised triangle groups. Two of the groups are non-uniform 
(and are discussed in \cite{MM2}) and the others are identified as fundamental 
groups of orbifolds obtained by Dehn surgeries on 2-bridge knots and links. 
As such they appear in \cite{HLM1} and our results establish a conjecture 
in that paper on the degree of the singular set of such orbifolds.

\medskip

As a basic reference to the deep relationships between arithmetic and 
hyperbolic geometry we refer to \cite{MR}.  We note here a few connections.  
In $\Isom^+\IH^2$,
 Takeuchi
\cite{Tak} identified all 82 arithmetic lattices generated by two elements
of finite order (equivalently arithmetic Fuchsian triangle groups), and, 
all arithmetic Fuchsian groups with two generators have been identified
(see  \cite{Tak2, MRos2}).  The
connections between arithmetic surfaces,  number theory and 
theoretical physics can be found in work of Sarnak and co-authors  eg.
\cite{Sar2,RS},  see also
\cite{katok}.  In
\cite{Vin1} Vinberg gave criteria for Coxeter groups in $\Isom\, \IH^n$ to be 
arithmetic. Such groups do not exist in the co-compact case for dimension $n \geq 30$,
\cite{Vin2}. It has now been established that there are finitely many 
conjugacy classes of maximal arithmetic Coxeter groups in all dimensions. 
There are two proofs see \cite{Agol} and \cite{Nik1} 
(previously established in two-dimensions in \cite{LMR}).

   Returning to dimension 3, the orientation-preserving subgroups of Coxeter
   groups for tetrahedra, some of which are generated by two elements of finite
order, which are arithmetic, are identified in \cite{Vin1} (see  \cite{MR2,
CM} for related results).
   Reid
\cite{R3} identified the figure eight knot complement as the only 
arithmetic knot complement.  The four  orientable hyperbolic 3-manifolds with 
fundamental group generated by a pair of parabolic elements which are 
arithmetic are two 
bridge knot and link complements \cite{GMM}.  The 14 finite co-volume Kleinian groups
with two generators, one of finite order, one parabolic, which are arithmetic 
are described in \cite{CMMO}.  An 
algorithmic approach to deciding if an orbifold obtained by 
$(p,0)-(q,0)$ surgery on a two bridge link (or knot) has arithmetic fundamental
groups was given in \cite{HLM1}.
 Arithmetic hyperbolic torus bundles are 
discussed in \cite{HLM,Bow}  and generalised triangle groups which are 
Kleinian in \cite{HMR,HLM2,Vin3}. Various extremal
groups have been identified as two-generator arithmetic; 
for instance the minimal 
volume non-compact hyperbolic 3--manifold and orbifold
\cite{Meyer, CaoM}.  
The Week's manifold is arithmetic, two-generator and conjecturally the 
minimal volume orientable hyperbolc 3-manifold \cite{CFJR}. 
The prime candidate for the minimal
volume orientable hyperbolic 3-orbifold is also arithmetic and generated
by two elements of finite order 
\cite{Chin2,GM3,GM4,MMarshall}.  

\medskip

Before precisely stating our main result let us say a few words about its 
proof  and why we have the restriction $p,q\geq 6$. In our work identifying 
the two-generator non-uniform lattices in \cite{MM2}, a key observation
was that the non-compactness hypothesis provided {\em a priori} knowledge 
that the underlying fields were quadratic  imaginary and the groups we were 
looking for were commensurable with Bianchi groups.  In this paper a major part of the work is to identify the underlying 
fields. For this we make use of some important
results of Stark \cite{Stark} and  Odlyzko \cite{Od},  
as well as results 
concerning the discriminants of number fields
of small degree such as those of Diaz Y Diaz and Olivier 
\cite{Di,CDO,DO}.  Using these bounds and some results from
the geometry of numbers and various discreteness criteria, we bound 
the degree of the fields in question and then, in
turn, bound the possible parameters for an arithmetic Kleinian group - once 
we have fixed the orders of the generators, the space of all discrete groups 
up to conjugacy is parameterised by a one complex-dimensional space.

Finally to identify all the groups we use  a
computer search to examine all algebraic integers in the field 
satisfying the given bounds and additional arithmetic restrictions on the real embeddings.  This procedure gives us a
relatively short list of candidate discrete groups which are now known to 
be subgroups of arithmetic Kleinian groups \cite{GMMR}.
We then use various ideas,  discussed in the body of the text, to 
decide if these groups are in fact arithmetic -  at issue here is the finiteness of the co-volume.

Thus we are able identify (up to
conjugacy) all the arithmetic Kleinian groups $\langle f,g
\rangle$ generated by an element
$f$ of order $p$  and $g$ of order $q$ with $p$ and $q$  at least 
$6$.  
Our results here also give the cases $p=2$,  $q\geq 6$ 
by the known result that a $(2,p)$-arithmetic hyperbolic lattice contains 
a $(p,p)$-arithmetic hyperbolic lattice with index at most two.   

  At present the remaining cases $p=2,3,4,5$ and $q\geq p$ 
seem computationally infeasible,  unless $q$ is
large enough - although we have  made some recent progress using the work of \cite{FR} on the most difficult case $p=2$ and $q=3$.  As the reader will come to realise,  the main problem here is finding effective bounds on the degree of the associated number fields.  
 
 \medskip  Here is our main result:

\begin{theorem}  Let  $\Gamma=\langle f,g \rangle$ be an arithmetic 
Kleinian group generated by elements of
order $p$ and $q$ with $p,q \geq 6$.  Then $p,q$ fall into one of the following 
4 cases:
\begin{enumerate}
\item $p=q=6$, \\ there are precisely  $12$ groups enumerated below in Table 1.
(See comments following the table).
\item $p=q=8$, \\ there is precisely one group obtained by (8,0) 
surgery on the knot 5/3.
\item $p=q=10$, \\there is precisely one group obtained by (10,0) 
surgery on each component of the link 13/5.
\item $p=q=12$, \\ there are precisely two groups obtained by (12,0) 
surgery on the knot 5/3  and on each component
of 8/3. 
\end{enumerate}
\end{theorem}
Here $r/s$ denotes the slope (or Schubert normal form) of a  two 
bridge knot or link,  \cite{BZ} \S 12.  Thus 5/3
denotes the well known figure eight knot complement.

\bigskip

  As a corollary we are able to verify a condition noticed by  Hilden, 
Lozano and Montesinos
\cite{HLM1} concerning the $(n,0)$ surgeries on two bridge link complements.
\begin{corollary}  Let $(r/s,n)$ denote the arithmetic hyperbolic orbifold 
whose underlying space is the 3-sphere and whose 
singular set is the 2 bridge knot or link with
slope
$r/s$ and has degree $n$.  Then
\begin{equation}
n \in \{2,3,4,5,6,8,10,12,\infty \}
\end{equation}
\end{corollary}

We also have the following,  slightly surprising, corollary.

\begin{corollary}\label{notrigrp}  There are no co-compact arithmetic generalised triangle 
groups with generators of orders at
least  $6$.
\end{corollary}

To prove the corollary we shall show later (see \ref{proofngtg}) that there 
cannot be another presentation of the same group on two generators of 
orders at least 6 as a generalised triangle group.

 \bigskip

In two dimensions,  there are in fact 32 arithmetic triangle groups 
with two generators of orders at least $6$ on
Takeuchi's list \cite{Tak}.  In three dimensions for each $p$ and $q$ 
($\min\{p,q\}\geq 3 $) there are infinitely many co-compact generalised 
triangle groups with a presentation of the form 
$\langle f,g: f^p=g^q=w(f,g)^2=1\rangle$ for certain words $w$ in $f$ and $g$.
Some of these are discussed in \cite{JR}.  Apparently as soon as 
$\min\{p,q\}\geq 6 $  none of these groups can be arithmetic.

\begin{table}
\centering
  \caption[]{Arithmetic groups with $p,q=6$ }
\begin{tabular}{|c|c|c|}
\hline
  $\gamma$ value & $k\Gamma$ & description of orbifold \\
\hline
   $ i\sqrt{3} $ &$z^2-z+1$ & $\Gamma_{21}$
\\
\hline
  $-1+i$ & $z^2+1$ & (6,0) surgery on 5/3 \\
\hline
  $-1$ & $z^2+z+1$ & $\Gamma_{20}$\\
\hline
   $1+3i$  &$z^2-2z+2$ & (6,0)-(6,0) surgery on link 24/7
\\
\hline
  $-1+i\sqrt{7}$ & $z^2-z+2$ & (6,0)-(6,0) surgery on link 30/11
\\
\hline
   $-2+i\sqrt{2}$  & $z^2+2$ & (6,0)-(6,0) surgery on link 12/5
\\
\hline
  $4.1096 - i\ 2.4317$ & $1+2z-3z^2+z^3$ &
(6,0) surgery on knot 65/51 \\
\hline
    $3.0674 -i\ 2.3277$ & $2-2z^2+z^3$ & (6,0) surgery on knot 13/3$^*$ \\
\hline
   $2.1244 -i\ 2.7466$ & $1+z-2z^2+z^3$ &  (6,0) surgery on knot 15/11 \\
\hline
    $1.0925 - i\ 2.052 $ & $1-z^2+z^3$ & (6,0) surgery on knot 7/3$^*$ \\
\hline
  $0.1240 -i\ 2.8365$ & $1+z-z^2+z^3$ & (6,0) surgery on knot 13/3$^*$\\
\hline
  $-0.8916 -i\ 1.9540 $ & $1+z+z^3$ & (6,0)-(6,0) surgery on link 8/5 \\
\hline
   $-1.8774 - i\ 0.7448$ & $1+2z+z^2+z^3$ & (6,0) surgery on knot 7/3$^*$ \\
\hline
    $-2.8846 -i\ 0.5897$ & $1+3z+z^2+z^3$ &
(6,0)-(6,0) surgery on link 20/9 \\
\hline
\end{tabular}
\end{table}
\bigskip

\noindent{\bf Notes:} Table 1,  and the results above,  were produced as follows.  The 
methods we outlined above and discuss in detail in 
the body of the paper produce for us all possible values of the trace 
of the commutator of a pair of
primitive elliptic generators of an arithmetic Kleinian group (the parameters) as well 
as an approximate volume for
the orbit space.   We then use Jeff Weeks' hyperbolic geometry 
package ``Snappea'' \cite{snap} to try and
identify the orbifold in question by surgering various two bridge 
knots and links and comparing
volumes.  Once we have a likely candidate,  we use the matrix presentation 
given by Snappea and verify that
the commutator traces are the same.  As these traces come as the 
roots of a monic polynomial with integer coefficients of modest
degree,  this comparison is exact.  Since this trace determines the 
group up to conjugacy,  we thereby identify the orbit space.
Conversely, once the two bridge knot or link and the relevant surgery 
is determined, a value of $\gamma$ can be recovered from the algorithm
in \cite{HLM1}.

Next,  if  $\alpha$ is a complex root of 
the given polynomial in Table 1,  then
$k\Gamma=\IQ(\alpha)$ and $\alpha(\alpha+1)=\gamma$.
Recall that the parameter $\gamma$ is determined by the Nielsen
equivalence class of a pair of generators of the groups.
In these tables an $*$ denotes that a Nielsen inequivalent pair of
generators  of order 6 (also listed in the table) gives rise to the 
same group.  Curiously these
examples were identified as follows.  Each is an index two subgroup of a group
obtained by $(2,0)$-surgery on one component $C_1$ and $(6,0)$-surgery on 
the other component $C_2$ of a
two bridge link,  in particular the links $7^{2}_{1}$ and $9^{2}_{1}$ in 
Rolfsen's tables \cite{Rolfsen}.  Suppose images of the meridians
are
$f$ of order six and
$g$ of order two.  Then $f$ and $gfg^{-1}$ are generators,  both of 
order six for the group in
question.  If however we do $(6,0)$-surgery on $C_1$ and $(2,0)$-surgery 
on $C_2$,  with images of
meridians being $f'$ and $g'$,  then $f'$ and $g'f'g'^{-1}$ give the 
same group,  but are not Nielsen
equivalent (as the $\gamma$ parameters are different).  Of course 
once identified,  one can use the
retriangluation procedure on Snappea to try to generate these different Nielsen
classes of generators (knowing they exist is a big incentive to 
retriangulate a few times).

The groups listed as $\Gamma_{20}$  and $\Gamma_{21}$ are the only 
non-compact examples and were
found in \cite{MM2}.  They have the following presentations
\[ \Gamma_{20} = \langle 
x,y:x^6=y^6=[x,y]^3=([x,y]x)^2=(y^{-1}[x,y])^2=(y^{-1}[x,y]x)^2=1\rangle 
\]
\[ \Gamma_{21} = \langle x,y:x^6=y^6=
(y^{-1}x)^2y[x^{-1},y][x,y][x,y^{-1}]x^{-1} = ([y^{-1},x]yx^2)^2 =1\rangle \]

\section{Two-generator Arithmetic Lattices}

The group of orientation preserving isometries of  the upper half-space 
model of $\IH^3$, 3-dimensional hyperbolic space,  is given by 
the group $\PSL(2, \IC)$, the 
natural action of its elements by linear fractional transformations 
on $\hat{\IC}$ 
extending to $\IH^3 = \IC \times \IR^+$ and preserving the metric of constant
negative curvature via the Poincar\'e extension.  

A subgroup $\G$ of
$\PSL(2,\IC)$ is said to be {\it reducible} if all elements have a common 
 fixed point in their action on $\hat{\IC}$ and $\G$ is otherwise {\it irreducible}. 
Also $\G$ is said to be {\it elementary} if it has a finite orbit in its action on
$\IH^3 \cup \hat{\IC}$ and $\G$ is otherwise {\it non-elementary}. 

\subsection{Parameters}
If $f \in \PSL(2,\IC)$ is 
represented by a matrix $A \in \SL(2,\IC)$ then the trace of $f$, $\tr(f)$,
is only defined up to a sign. However, if $[f,g]=fgf^{-1}g^{-1}$ denotes the commutator of $f$
 and $g$, then $\tr[f,g]$ is well-defined and, furthermore,
the two generator group $\langle f,g\rangle $ is reducible if and only if
$\tr[f,g]=2$. For a two-generator group $\langle f,g\rangle $ the three complex numbers
 $(\gamma(f,g),\beta(f), \beta(g))$   
 \begin{equation}\label{betadef} \beta(f) = \tr^2(f) - 4, \hskip10pt \beta(g) = \tr^2(g) - 4, \hskip10pt \gamma(f,g) = \tr[f,g]-2\end{equation}
 are well-defined by $f,g$ and form the 
{\it parameters} of the group $\langle f,g\rangle $. They define $\langle f,g\rangle $ uniquely up to 
conjugacy provided $\langle f,g\rangle $ is irreducible, that is $\gamma(f,g) \neq 0$, see \cite{GM0}.

Now suppose that $f$ and $g$ have finite orders $p$ and $q$ respectively
where we can assume that $p \geq q$. In considering the group 
$\G = \langle  f , g \rangle $ we can assume that $f$ and $g$ are primitive elements and so $\G$ has parameters
\begin{equation}   (\gamma, -4 \sin^2 \pi/p, -4 \sin^2 \pi/q).
\end{equation}
(Where there is no danger of confusion, we will abbreviate $\gamma(f,g)$ 
simply to $\gamma$.)  For fixed $p,q$, any $\gamma \in \IC\setminus \{ 0\}$
uniquely determines the conjugacy class of such a group $\G = \langle  f , g \rangle $. We say
$\G$ is {\it Kleinian} if it is a discrete non-elementary subgroup of $\PSL(2,\IC)$. For fixed $p$ and $q$ it is an elementary consequence of a theorem of J\o rgensen \cite{Jorg} that the 
set of all such $\gamma$ is closed and computer generated pictures suggest 
that it is highly fractal in nature - for instance the Riley slice, 
corresponding to two parabolic generators, would correspond to $p=q=\infty$. 

The cases where $\gamma$ is real have been investigated
in \cite{Klim,KlimKop,MM3}
 We have shown in \cite{MM}, that for each pair 
$(p,q)$ there are only finitely many $\gamma$ in $\IC$ 
 which yield  arithmetic
Kleinian groups and for all but a finite number of pairs $(p,q)$, that 
finite number is zero. It is our aim here to determine all $\gamma$
such that $\G$ is an arithmetic Kleinian group (i.e. a 3-dimensional arithmetic hyperbolic lattice) with $p,q \geq 6$ and to obtain a geometric description of these groups.

\subsection{Arithmetic Kleinian Groups}

For detailed information on arithmetic Kleinian groups see \cite{Bo,Vig,MR}.  For completeness,  and since we will rely heavily on these results, we recall here some basic facts.

Let $k$ be a number field and for each place $\nu$ of $k$, let $k_{\nu}$
denote the completion of $k$ with respect to the metric on $k$ induced by the 
valuation $\nu$. For each Galois monomorphism $\sigma : k \rightarrow \IC$, there
is an Archimedean valuation given by $| \sigma(x) |$ and if $\sigma(k) \subset \IR$ then $k_{\nu} \cong \IR$ and, if not, each complex conjugate
 pair forms a place
and $k_{\nu} \cong \IC$. The other valuations are ${\cal P}-$adic and 
correspond to prime ideals ${\cal P}$ of $R_k$. The fields $k_{\nu} = 
k_{\cal P}$ are finite extensions of the $p$-adic numbers $\IQ_p$. 
Let $A$ be a quaternion algebra over $k$ and let $A_{\nu} = A \otimes_k
k_{\nu}$ so that $A_{\nu}$ is a quaternion algebra over the local
field $k_{\nu}$. For $k_{\nu} \cong \IC$, then $A_{\nu} \cong M_2(\IC)$
but for all other places there are just two quaternion algebras over  
each local field one of which is $M_2(k_{\nu})$ and the other is a unique
quaternion division algebra over $k_{\nu}$. 

 We say that $A$ is {\it ramified} at $\nu$
if $A_{\nu}$ is a division algebra. The set of places at which $A$ is ramified 
is finite of even cardinality   and is called the    
{\it ramification set} of $A$, denoted by $\Ram(A)$.  The ramification set determines the isomorphism class 
of $A$ over $k$. We also denote the set of Archimedean ramified places by
$\Ram_{\infty}(A)$ and the non-Archimedean or finite places at which $A$
is ramified by $\Ram_f(A)$. Now as a quaternion algebra $A$ has a basis of the form $1,i,j,ij$
where $i^2 = a, j^2 = b$ and $ij = -ji$, with $a,b \in k^*$. It can thus be
represented by a {\it Hilbert symbol} $\ds{\HS{a}{b}{k}}$. If the real
place $\nu$ corresponds to the embedding $\sigma : k \rightarrow \IR$,
 then 
 \[ A_{\nu} \cong 
\ds{\HS{a}{b}{k} \otimes_k k_{\nu} \cong \HS{\sigma(a)}{\sigma(b)}{\IR}} \]
and $A$ will be ramified at $\nu$ if and only if $A_{\nu}$ is isomorphic to
Hamilton's quaternions. This occurs precisely when both $\sigma(a)$ and 
$\sigma(b)$ are negative.

Now assume that $k$ has exactly one complex place and that the quaternion algebra
$A$ is ramified at least at all the real places of $k$. Let ${\cal O}$ be an order in $A$ and let ${\cal O}^1$ denote the elements of norm 1. In these 
circumstances there is a $k$-embedding $\rho : A \rightarrow M_2(\IC)$
and the group $P\rho({\cal O}^1)$ is a Kleinian group of finite co-volume.
The set of {\it arithmetic Kleinian groups} is the set of Kleinian groups which are
commensurable with some such $P\rho({\cal O}^1)$. 

It is  our aim to 
identify all (conjugacy classes of) arithmetic Kleinian groups generated by two elements 
of finite order.  To do this we use the identification theorem below
which gives a method of identifying
arithmetic Kleinian groups from the elements of the given group. 

We require the following preliminaries. Let
$\G$ be any non-elementary finitely-generated subgroup  of $\PSL(2,\IC)$. 
Let $\G^{(2)} = \langle  g^2 \mid g \in \G \rangle $ so that $\G^{(2)}$ is a subgroup of finite
index in $\G$. Define
\begin{equation} \left.
\begin{array}{lll}
k\G & = & \IQ(\{ \tr(h) \mid h \in \G^{(2)} \}) \\ 
A\G  & = & \{ \sum a_i h_i \mid a_i \in k\G, h_i \in \G^{(2)} \}
\end{array}\;\;\;\;  \right\}
\end{equation}
where, with the usual abuse of notation, we regard elements of $\G$ as matrices,
so that $A\Gamma \subset M_2(\IC)$. 

Then $A\G$ is a quaternion algebra over $k\G$ and the pair $(k\G, A\G)$ is 
an invariant of the commensurability class of $\G$. If, in addition, $\G$ 
is a Kleinian group of finite co-volume then $k\G$ is a number field.

We state  the identification  
theorem as follows:
\begin{theorem} \label{idthm} Let $\G$ be a  subgroup of 
$\PSL(2,\IC)$ which is finitely-generated and non-elementary. Then $\G$
is an arithmetic Kleinian group if and only if the following conditions all hold:
\begin{enumerate}
\item  $k\G$ is a number field with exactly one complex place,
\item for every $g \in \G$, $\tr(g)$ is an algebraic integer,
\item $A\G$ is ramified at all real places of $k\G$.
\item $\G$ has finite co-volume.
\end{enumerate}
\end{theorem}
It should be noted  that the first three conditions together imply that $\G$ is 
Kleinian, and without the fourth condition, are  sufficient to imply that $\G$ is a subgroup of an arithmetic Kleinian group.

The first two conditions clearly depend on the traces of the elements of $\G$.
In addition, we may also find a Hilbert symbol for $A\G$ in terms of the 
traces of elements of $\G$ so that the third condition also depends on the traces 
(for all this, see \cite{MR1},\cite[Chap. 8]{MR}).

\subsection{Two-generator arithmetic groups}

We now  suppose that $\G$ is generated by two elements $f,g$ of orders $p$ and $q$ 
respectively where $p \geq q$.  We have noted that the conjugacy
class of $\G$ is uniquely determined by
the single complex parameter $\gamma$. We now show how the first three conditions of
Theorem \ref{idthm} can be equivalently expressed in terms of $\gamma$.  This
is not true of the fourth condition, but for $\G$ to have finite co-volume places
some necessary conditions on $\gamma$ (see \S 3 below).

Note that $\tr f = \pm 2\cos(\pi/p)$ and $\tr g= \pm 2\cos(\pi/q)$ are algebraic integers and recall that the traces
of all elements in $\langle  f , g \rangle $ are integer polynomials in $\tr f, \tr g$
and $\tr fg$. Now the Fricke identity states
\begin{equation}\label{Fricke}
\gamma=\gamma(f,g) = \tr^2 f + \tr^2 g + \tr^2 fg - \tr f \tr g \tr fg -4.
\end{equation}
Thus  $\tr fg$ is an algebraic integer if and only if $\gamma$ is an algebraic 
integer so that the second condition of Theorem \ref{idthm} is equivalent, in these two-generator cases, 
to requiring that $\gamma$ be an algebraic integer.

Now suppose that $p,q \geq 3$. Throughout, we denote $\beta(f), \beta(g)$ (see (\ref{betadef})) by 
$\beta_1, \beta_2$ respectively so that 
\[ \beta_1 = -4 \sin^2 \frac{\pi}{p},\hskip10pt  \beta_2 
= -4 \sin^2 \frac{\pi}{q}, \hskip10pt  \beta_1+4 = 4 \cos^2 \frac{\pi}{p}, \hskip10pt  \beta_2 + 4 = 4 \cos^2 \frac{\pi}{q} \] 
Now $k\G = \IQ(\tr^2 f, \tr^2 g, \tr f \tr g \tr fg)$ 
(see for instance \cite[Chap.3]{MR}). We consistently use $L$ to denote 
the totally real subfield
\[ L = \IQ( \tr^2 f, \tr^2 g) = \IQ(\beta_1, \beta_2) = 
\IQ(\cos \frac{2\pi}{p}, \cos \frac{2\pi}{q}) \]
Thus $k\G = L(\lambda)$ where $\lambda = \tr f \tr g \tr fg$. From the Fricke 
identity
(\ref{Fricke}) and $\tr^2(fg)=\lambda^2/(\beta_1+4)(\beta_2+4)$ 
we deduce that $\lambda$ satisfies the quadratic equation
\begin{equation}\label{eqn5}
x^2 -  (4+\beta_1)(4+\beta_2) \, x 
  + (4+\beta_1)(4+\beta_2)(\beta_1 + \beta_2 + 4 -\gamma) = 0,
\end{equation}
and that $[k\G : L(\gamma)]\leq 2$.

\bigskip\noindent\framebox[1.2\width]{{\bf $\gamma(f,g)\in\IR$} }

\bigskip

\noindent  Let us at this point remove the inconvenient case that $\gamma(f,g)$ is real as this case complicates our discussion. 
Suppose then that $\gamma\in \IR$. In the next section (see (\ref{eqn21})), it will
be shown that, for $\G$ to have finite co--volume we must have
\[ -4 < \gamma < 4(\cos \pi/p +   \cos \pi/q)^2 \]
 Now if $\gamma \geq 0 $, then for any Kleinian group
$\G = \langle f, g\rangle $ with $o(f) = p, o(g)=q$ and $\gamma(f,g) = \gamma$,  
$\G$ has an invariant plane \cite{Klim,MM3} and so, as the reader can easily verify,  
cannot have finite co-volume and hence cannot be an arithmetic Kleinian group. Thus $-2 < \tr[f,g] < 2$ 
and so, whenever $\G$ is 
discrete and finite co-volume the commutator $[f,g]$ must be 
elliptic. All such groups, arithmetic or otherwise, have been determined 
in \cite{MM3}.  There are precisely nine such groups which 
are arithmetic, all have $p,q \leq 6$ and there is only one with $p=q=6$.

\medskip

{\em Thus we assume henceforth that $\gamma$ is not real. }

\medskip
 
 Then $k\G$ will be  a number field with one complex place if and only if $L(\gamma)$
has one complex place and the quadratic at (\ref{eqn5}) splits into linear factors over
$L(\gamma)$. This implies that, if $\tau$ is any real embedding of 
$L(\gamma)$, then the image of the discriminant of  (\ref{eqn5}), which is 
$(4+\beta_1)(4+\beta_2)(\beta_1 \beta_2 + 4 \gamma)$, 
under $\tau$ must be positive. Clearly this is equivalent to requiring that
\begin{equation}\label{eqn6}
\tau( \beta_1 \beta_2 + 4 \gamma) > 0.
\end{equation}
Thus $k\G$ has one complex place if and only if (i) $\IQ(\gamma)$ 
has one complex 
place, (ii) $L \subset \IQ(\gamma)$, (iii) for all real embeddings $\tau$ of $\IQ(\gamma)$,
(\ref{eqn6}) holds and (iv) the quadratic at (\ref{eqn5}) factorises over $\IQ(\gamma)$.

 Now, still in the cases where $p,q > 2$,(\cite[\S 3.6]{MR})  
\begin{equation}\label{eqn7}
A\G = \HS{\beta_1(\beta_1+4)}{ (\beta_1+4)(\beta_2+4)\,\gamma}{k\G}.
\end{equation}
Under all the real embeddings of $k\G$, the term  $\beta_1(\beta_1+4)$ is 
negative and 
$(\beta_1+4)(\beta_2+4)$ is positive. Thus $A\G$ is ramified at
all real places of $k\G$ if and only if, under any real embedding $\tau$
of $k\G$,
\begin{equation}\label{eqn8}
\tau(\gamma) < 0.
\end{equation}
Thus, summarising, we have the following theorem which we will use to determine the possible $\gamma$ values for the groups we seek.
\begin{theorem}\label{2genthm}
Let $\G = \langle  f , g \rangle $ be a non-elementary
subgroup of $\PSL(2,\IC)$ with $f$ of order $p$ and $g$ of order $q$, $p \geq q \geq 3$.  Let $\gamma(f,g) = \gamma \in \IC \setminus \IR$. Then $\G$ is an arithmetic Kleinian group if and only
if 
\begin{enumerate}
\item $\gamma$ is an algebraic integer,
\item $\IQ(\gamma) \supset L = \IQ(\cos 2 \pi/p, \cos 2 \pi/q)$ and $\IQ(\gamma)$ is a number field with exactly one complex place,
\item if $\tau : \IQ(\gamma) \rightarrow \IR$ such that $\tau |_L = \sigma$, then
\begin{equation}\label{eqn10}
 - \sigma( \frac{\beta_1 \beta_2}{4}) < \tau(\gamma) < 0,
\end{equation}
\item the quadratic polynomial at (\ref{eqn5})   factorises over $\IQ(\gamma)$,
\item $\G$ has finite co-volume.
\end{enumerate}
\end{theorem}

Any non-elementary subgroup $\G = \langle  f, g \rangle $ of $\PSL(2,\IC)$ where $o(f) = 
o(g) = p > 2$ is contained as a subgroup of index at most 2 in a group 
$\G^* = \langle  h, f \rangle $ where $o(h) = 2$ with 
\begin{equation} \label{eqn11}
\gamma(f,g) = \gamma(h,f) ( \gamma(h,f)  - \beta_1)
\end{equation}
and conversely (see \cite{GM1}). Thus, $(k\G, A\G) = (k\G^*, A\G^*)$, since these are 
commensurability invariants, and so we can obtain necessary and sufficient 
conditions for arithmeticity of $\G$ in terms of $\gamma = \gamma(h,f)$
where $o(h)=2, o(f)=p > 2$. In this case, $k\G^* = \IQ( \tr^2 f, \gamma)
= L(\gamma)$ (see \cite{MR}) and 
\begin{equation}
A\G^* = \HS{\beta_1(\beta_1+4)}{\gamma\,(\gamma - \beta_1)}{k\G^*}.
\end{equation}
Arguing as above, we have 
\begin{theorem}\label{2gen*}
Let $\G^* = \langle  h,f \rangle $ be a non-elementary
subgroup of $\PSL(2,\IC)$ with $h$ of order $2$ and $f$ of order $p >2$.  Let  $\gamma(h,f) = \gamma \in \IC \setminus \IR$.
Then $\G^*$ is an arithmetic Kleinian group if and only if
\begin{enumerate}
\item $\gamma$ is an algebraic integer,
\item $\IQ(\gamma) \supset L = \IQ(\cos 2 \pi/p)$ and $\IQ(\gamma)$ is 
a number field with exactly one complex place,
\item if $\tau : \IQ(\gamma) \rightarrow \IR$ such that $\tau|_L = \sigma$ then
\begin{equation}\label{eqn13}
 \sigma(\beta_1) < \tau(\gamma) < 0,
\end{equation}
\item $\G^*$ has finite co-volume.
\end{enumerate}
\end{theorem} 

\medskip
Implementation of the fourth condition of Theorem \ref{2genthm} can be 
simplified as follows: 
Suppose that $m(x)$, the minimum polynomial of $\gamma$ over $L$, has the form
$x^r + a_{r-1} x^{r-1} + ... + a_0.$   From  our usual expression for 
$\gamma$  at (5), we have: 
\begin{equation}\label{eqn14}
 \tr^2 f \tr^2 g \;  \gamma = (\tr f  \tr g  \tr fg)^2 -  \tr^2 f  \tr^2 g (\tr f \tr g  \tr fg) +  \tr^2 f  \tr^2 g ( \tr^2 f + \tr^2 g  - 4).
 \end{equation}
 That is
$ b  \gamma     = \lambda^2 -  b  \lambda  + c $
where  $b$ and $c$ are integers in $L$.  Next,  substituting in $m(x)$ and 
clearing denominators gives $(\lambda^2 - b \lambda + c)^r  +  a_{r-1} b (\lambda^2 - b \lambda +c)^{r-1} + ..... + a_0  b^r  = 0 $
which is a monic polynomial in $\lambda$ of degree $2r$ with 
coefficients integers in $L$.  We define the polynomial
\[ M(y) =(y^2 - by + c)^r  +  a_{r-1} b (y^2 - b y +c)^{r-1} + ..... + a_0  b^r  \]
simply replacing $\lambda$ by $y$. 

Since  $\IQ(\lambda) = \IQ(\gamma)$,  then $\lambda$ is an algebraic integer 
in $k\G$ which has a minimum polynomial over $L $ which is monic with 
integer coefficients. This must also be true of the ``other" root  
$\lambda' = b - \lambda$. So the two factors of the polynomial $M(y)$ 
have coefficients which are integers in $L.$  Hence the fourth condition of 
Theorem \ref{2genthm} is equivalent to 
\begin{lemma}\label{condition4'} The polynomial   $M(y)$ factors over $L$ into two monic  factors both
of degree $r$ and having integral coefficients (in $L$)
\end{lemma}
A slight simplification of this occurs in the cases where $(p,q)>2$. In these cases,
$ \dis{a = 8 \cos \frac{\pi}{p} \, \cos \frac{\pi}{q} \, \cos ( \frac{\pi}{p} + \frac{\pi}{q})}$ is an algebraic integer in $L$. If we set $\epsilon = 
\lambda - a$, then $\IQ(\lambda) = \IQ(\epsilon)$ and equation (14) 
takes the form 
\begin{equation}\label{eqn15}
   b \gamma = \epsilon ( \epsilon - c)
   \end{equation}
   where $b= 16 \cos^2 \pi/p \cos^2 \pi/q$, $c= 4 \sin 2 \pi/p \sin 2 \pi/q$
   are integers in $L$. We can use this factorisation in $m(x)$ to obtain the corresponding result to Lemma 2.4.

See (21) later for an example of this condition applied - using an integral basis for $L$ this condition can be rewritten to assert the existence of a solution in rational integers of a nonlinear system of equations.  Since our methods are to deduce the possible minimum polynomials of $\gamma$ over $L$, this alternative formulation can be readily computationally implemented.
Note that when $p=q$, 
\[ \epsilon = - 4 \cos^2\frac{\pi}{p} \,\, \gamma(h,f) \]  
and  (\ref{eqn11}) is a special case of (\ref{eqn15}) and hence of (\ref{eqn14}).

\section{Free Products}

As we have noted, the first four conditions of  Theorem \ref{2genthm}  on $\gamma$ are  sufficient to imply that $\G$ is a subgroup of an arithmetic Kleinian group.  However many of the 
groups satisfying these four conditions will  be isomorphic to the free 
product $\langle  f \rangle  * \langle  g \rangle $ and so cannot be arithmetic Kleinian groups
as they must fail to have finite co-volume. To eliminate these groups we now seek conditions on 
$\gamma$ which force a discrete group $\G = \langle  f , g \rangle $ to be a free product.
Moreover, we will  extend the methods of \cite{MM} to enumerate
the parameters $\gamma$ which give rise to arithmetic Kleinian groups
by obtaining bounds which involve the discriminant of the power basis
of $\IQ(\gamma)$ over $L$ determined by $\gamma$. For this purpose, and also
for other methods to be used in the enumeration, we want to obtain as stringent
bounds as possible on $| \gamma |, \Im(\gamma), \Re(\gamma)$.
The extreme values of these are attained within a contour $\Omega_{p,q}$ in the
$\gamma$-plane. We thus obtain bounds which are simple functions of one 
variable which, for each pair $(p,q)$ can be (computationally) maximised.
\bigskip

Define 
\[ A = \begin{pmatrix}  \cos \pi/p & i \sin \pi/p \\ i \sin \pi/p & 
\cos \pi/p \end{pmatrix}, \quad B = \begin{pmatrix}
 \cos \pi/q & iw \sin\ \pi/q \\ i w^{-1} \sin \pi/q & \cos \pi/q \end{pmatrix}. \]
Then if $\G = \langle  f,g\rangle $ is a non-elementary Kleinian group with $o(f)=p,
o(g) = q$, where $p \geq q \geq 3$ then $\G$ can be normalised so that 
$f,g$ are represented by the matrices $A,B$ respectively. The parameter
$\gamma$ is related to $w$ by
\begin{equation}\label{eqn16}
\gamma = \sin^2 \frac{\pi}{p} \, \sin^2 \frac{\pi}{q} \; (w - \frac{1}{w})^2.
\end{equation}
Given $\gamma$, we can further normalise and choose $w$ such that $| w | 
\leq 1$ and ${\rm Re}(w) \geq 0$. 

It is convenient here to also consider the cases where $\G^* = \langle h,f \rangle$
with $o(h)=2, o(f)=p$ as discussed in Theorem \ref{2gen*} so that in this section we will 
allow  $q$ to be equal to 2. 

We recall the {\em isometric circles} of a linear fractional transformation
 \[ g(z)=\frac{az+b}{cz+d} \approx  \begin{pmatrix}  a &b \\ c & 
d \end{pmatrix} \in PSL(2,\IC), \hskip10pt c\neq 0 \]
are the pair of circles
\[ I(g) = \{z : |cz+d|=1\},\hskip20pt  I(g^{-1}) = \{z : |cz-a|=1\} \]
Notice that $I(g)=\{|g'(z)|=1\}$ and $I(g^{-1})=\{|(g^{-1})'(z)|=1\}$ and that $g$ maps the exterior of $I(g)$ to the interior of $I(g^{-1})$.

The Klein combination theorem, (see \cite{Mas} for this and important 
generalisations) can be used to establish the following well known fact:  If the isometric circles of $g$  lie inside the 
intersection of the disks bounded by the isometric circles of $f$, then 
$\langle  f,g \rangle  \cong \langle f\rangle  * \langle  g \rangle $.
(See the illustrative examples in Diagram 1, where this situation holds
in case 1 but not in case 2.)

\scalebox{0.75}{\includegraphics*[viewport=30 200 500 440]{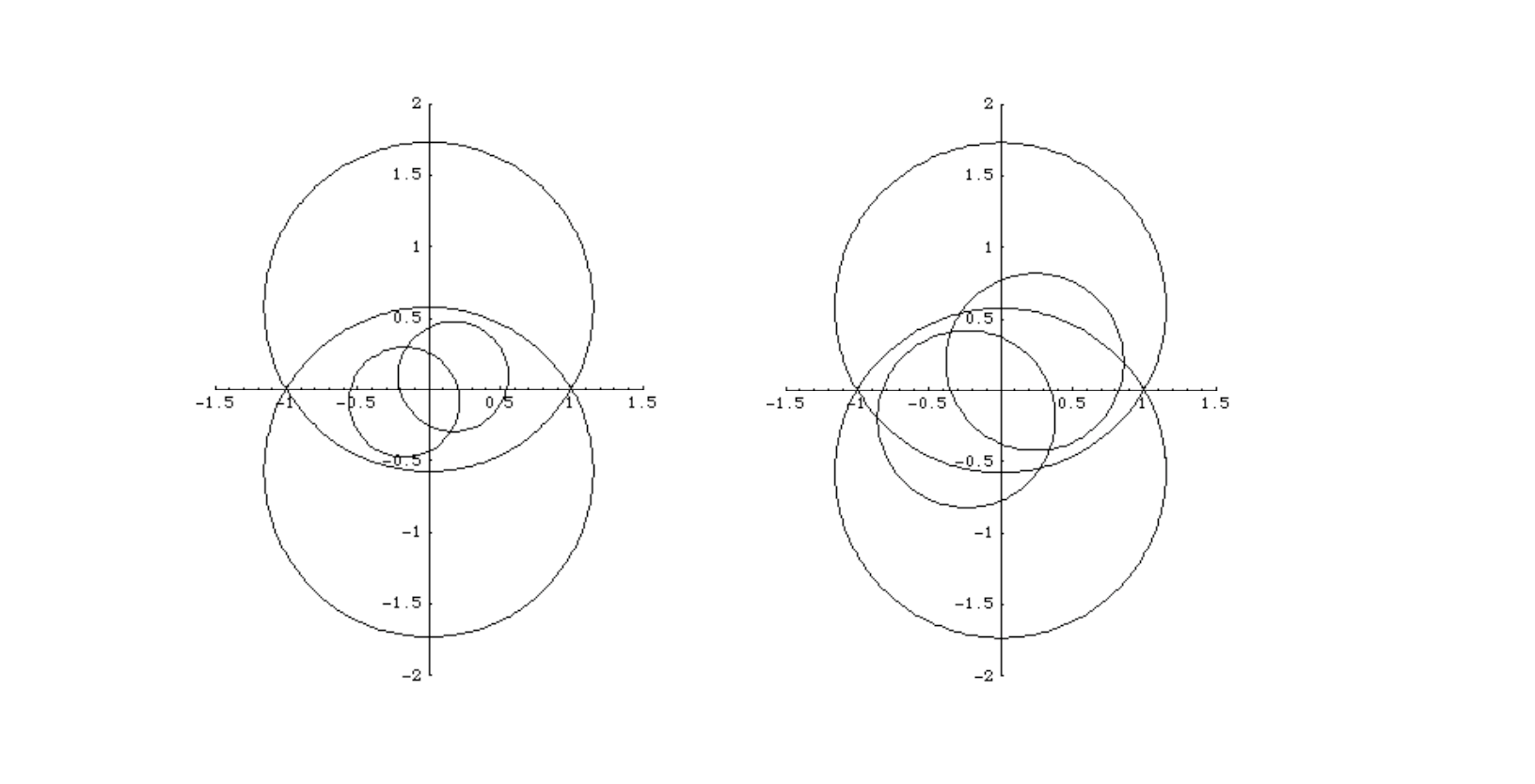}}

 \begin{center}  {\bf Diagram 1. $p=q=3$ isometric circles};\\ 1. non-intersecting ($\gamma=-4 + 4 i$) 2. intersecting ($\gamma=-1.5 + 1.75 i$) 
 \end{center}

This geometric configuration occurs precisely when
\begin{equation}\label{eqn17}
| i w \cot \pi/q + i \cot \pi/p | + \frac{| w |}{\sin \pi/q} \leq \frac{1}{\sin \pi/p}.
\end{equation}
As $w$ traverses the boundary of the region described by (\ref{eqn17}), 
then $\gamma$ traverses a contour $\Omega_{p,q}$, so that, when $\gamma$
lies outside this, the corresponding group will be a free product. The 
general shape of such a contour is illustrated by the case exhibited
in Diagram 2.
\begin{center}
\includegraphics*[viewport=-20 230 350 400]{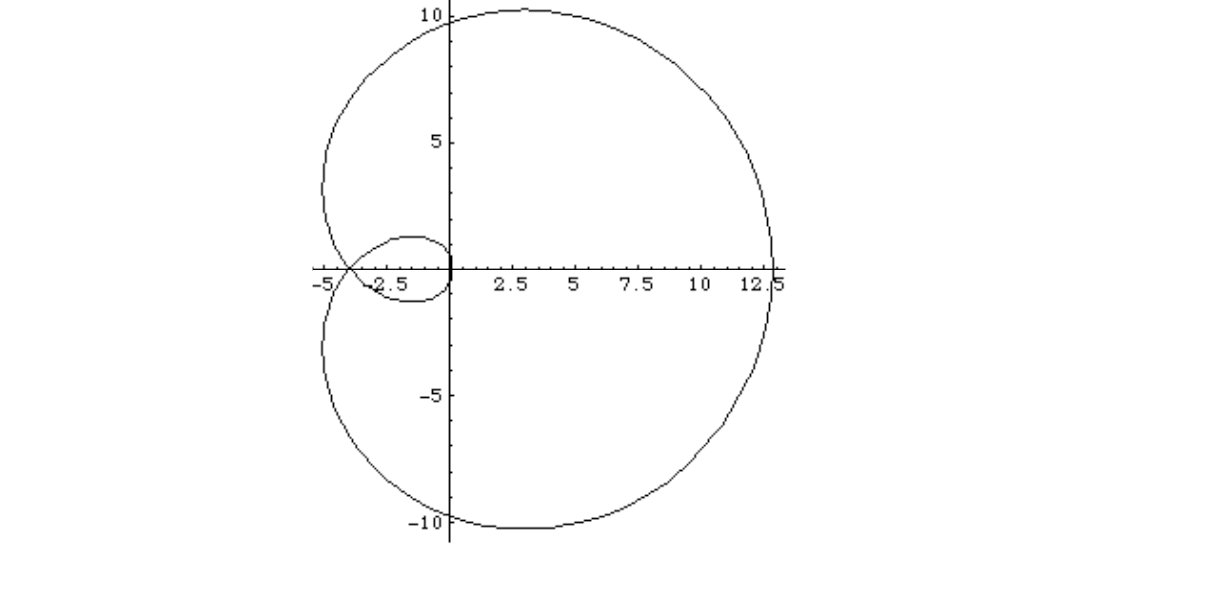}\newline
{\bf Diagram 2}\quad  $\Omega_{8,6}$
\end{center}
More specifically, let 
$w = r {\rm e}^{i \theta}$, and define  
$c(p,q) = \cos \pi/p \, \cos \pi/q$,
$s(p,q) = \sin \pi/p \, \sin \pi/q$. Then on the boundary of the region 
defined by (\ref{eqn17}),
\begin{equation}\label{rexprn} r^2 + \frac{1}{r^2} = \frac{4(1+c(p,q)\cos \theta)^2}{s(p,q)^2} - 2 \end{equation}
Since
$$
\gamma  =    
s(p,q)^2\,[(r^2 + r^{-2})\cos 2 \theta - 2 + 
(r^2 - r^{-2}) i \sin 2 \theta],
$$
and we can assume that $|w|<1$ and $\Re(w) \geq0$, we set $\cos \theta = t 
\in [0,1]$ and obtain
\begin{eqnarray}\label{Omega}
\Omega_{p,q}(t)= & 4(2t^2-1)(1+t c(p,q))^2-4t^2 s(p,q)^2  \\ 
  & -8t \sqrt{1-t^2} (1+t c(p,q)( \sqrt{(1+t c(p,q))^2-s(p,q)^2} i \nonumber 
  \end{eqnarray}
  It is clear that the real part of $\gamma$ takes its maximum value for $t=1$
  and so 
\begin{equation}\label{eqn20}
 \Re(\gamma) \leq 4( \cos \pi/p + \cos \pi/q)^2.
\end{equation}
and for each $(p,q)$ its minimum value can be computed from this formula.
Note that, if $\gamma$ is real, then
\begin{equation}\label{eqn21}
-4 \leq \gamma \leq 4(\cos \pi/p + \cos \pi/q)^2.
\end{equation}
which gives us the estimate we used earlier to handle the case $\gamma\in\IR$. 

More generally, for $\gamma \in \Omega_{p,q}(t)$, we have
$$
| \gamma | = 4[(1 + t c(p,q))^2  - t^2 s(p,q)^2].
$$
When  $c(p,q) \geq s(p,q)$, which occurs in particular when $p,q \geq 6$,
\begin{equation}\label{eqn22}
| \gamma | \leq 4( \cos \pi/p + \cos \pi/q)^2.
\end{equation}
Also, in the case $(p,2)$, we have $| \gamma | \leq 4$.
 Finally, note that 
$$| \gamma + 4 s(p,q)^2 | + | \gamma| = 2 s(p,q)^2 ( r^2 + r^{-2}).$$
From the expression for $r^2+r^{-2}$ at (\ref{rexprn}) above, this clearly takes its maximal value
when $\cos \theta = 1$. Thus if $x$ is a real number in the interval
$[-\beta_1 \beta_2 / 4,0] = [-4 s(p,q)^2, 0]$, then 
\begin{equation}\label{eqn23}
| \gamma - x | \leq 4(1 + \cos \pi/p \cos \pi/q)^2.
\end{equation}

\section{The possible values of $(p,q)$}

From Theorem \ref{2genthm}, we note, first, that $\gamma$ is an algebraic integer,
secondly, that $\IQ(\gamma)$ has exactly one complex place and thirdly, 
that $\IQ(\gamma)$  must contain $L = \IQ(\cos 2 \pi/p,
\cos 2 \pi/q)$. Let $[\IQ(\gamma) : L ] = r$. We now make use of these facts,
together with  the inequalities that $\gamma$ and its conjugates must satisfy
given in \S 2 and \S 3  to produce a list of possible
values for the triple $(p,q,r)$ for which there may exist a $\gamma$-parameter 
corresponding to an arithmetic Kleinian group which is not obviously free 
using the criteria from \S 3.  Further, if $(p,q,r)$ 
does not appear on this list there cannot be any corresponding arithmetic
Kleinian groups (see Table 5). The list obtained in Table 5 is produced
by refining a basic list in \S 4.1 using arguments on the norm and discriminant,
each stage being implemented by an elementary program in Maple. The 
finiteness of such a list was established in 
\cite{MM} and the starting point here uses the crude estimate obtained
in \cite{MM} that
 $p,q \leq 120$. In producing our lists, we assume that $p,q \geq 6$ 
although the methods apply for $p,q \geq 3$.

\subsection{Norm method} 

Let $N$ denote the absolute norm $N : \IQ(\gamma) \rightarrow \IQ$ 
and, as before, $L = \IQ(\cos 2 \pi/p, \cos 2 \pi/q)$.
If $(p,q)>2$, then $L = \IQ(\cos 2 \pi/M)$ where $M$ is the least common 
multiple of $p$ and $q$  and otherwise $L$ is of index 2 in that field.
Thus if $\mu = [L : \IQ]$, then
$$ \mu = \left\{ \begin{array}{ll}
             \phi(M)/2 & {\rm if~}(p,q)>2 \\
             \phi(M)/4 & {\rm if~}(p,q)\, |\, 2.
           \end{array}
         \right. $$
 The field $L$ is totally
real and the embeddings $\sigma : L \rightarrow \IR$ are defined by 
$$\sigma( \cos \frac{2 \pi}{p})  =   \cos \frac{2 \pi j}{p}, 
\sigma(\cos \frac{2 \pi}{q}) = \cos \frac{2 \pi j}{q}~{\rm where~}(j,pq)=1.$$
 Let us denote
these embeddings by $\sigma_1, \sigma_2, \ldots , \sigma_{\mu}$, 
with $\sigma_1 = {\rm Id}$.
Since  $\gamma$ is an algebraic integer
$| N(\gamma)| \geq 1$ and $N(\gamma) = \gamma \bar{\gamma} \prod_{\tau} \tau(\gamma)$ where $\tau$ runs over the $r \mu - 2 $ real embeddings of 
$\IQ(\gamma)$. 
If $\tau |_L = \sigma_i$, then by (\ref{eqn10}) 
$$ - \frac{\sigma_i(\beta_1 \beta_2)}{4} < \tau(\gamma) < 0$$
and from (\ref{eqn22}), $| \gamma| < 4(\cos \pi/p + \cos \pi/q)^2$. Thus we obtain
\begin{equation}{\label{Norm1}}
1\leq | N(\gamma)| \leq  
    16(\cos \pi/p + \cos \pi/q)^4 \, (4 s(p,q)^2)^{-2} \prod_{j=1}^{\mu} \left( \frac{\sigma_j(\beta_1 \beta_2)}{4} \right)^r.
\end{equation}
Now letting 
\[ \delta_n = \left\{ \begin{array}{ll}
                1 & {\rm if~}n \neq p^{\alpha}, \;\; p~{\rm a~prime} \\
               p & {\rm if~}n = p^{\alpha}, \;\; p ~{\rm a~prime},
                 \end{array}
        \right.  \]
then, (see \cite{MM}), if $\delta_{n,m} = \delta_n^{2/\phi(n)}\, \delta_m^{2/\phi(m)}$, 
\begin{equation}
\prod_{j=1}^{\mu}  \sigma_j(\beta_1 \beta_2) = 
\delta_{p,q}^{\mu}.
\end{equation}
Thus, taking logs, for a triple $(p,q,r)$ to give rise to a $\gamma$ which  
 represents an arithmetic Kleinian group it must satisfy the inequality
\begin{equation}\label{eqn30}
r \mu \leq 4 \log \left[ 
\frac{\cos \pi/p + \cos \pi/q}{\sin \pi/p \sin \pi/q} \right] /\log (4/\delta_{p,q})\end{equation}
Note that $r \geq 2$ and $6 \leq p,q \leq 120$ so that we can determine the triples 
for which (\ref{eqn30}) holds by obtaining the values of $p$ and $q$ and an upper 
bound for the related value of $r$. This produces a list of 
86   entries shown in Table 2,
which, for future reference, we call 
the {\it Norm List}.  
\begin{table}
\centering
\begin{tabular}{||c|c|c||c|c|c||c|c|c||c|c|c||}
\hline 
p&q&r&p&q&r&p&q&r&p&q&r \\ \hline
6&6&5&7&6&3&7&7&33&8&6&4 \\ \hline
8&7&4&8&8&7&9&6&3&9&7&3 \\ \hline
9&8&2&9&9&5&10&6&3&10&7&2 \\ \hline
10&8&2&10&10&4&11&6&2&11&7&2 \\ \hline
11&8&2&11&11&5&12&6&3&12&7&2 \\ \hline
12&8&2&12&9&2&12&10&2&12&12&4 \\ \hline
13&7&2&13&13&4&14&6&2&14&7&5 \\ \hline
14&14&3&15&6&2&15&10&2&15&15&2 \\ \hline
16&6&2&16&8&3&16&16&3&17&17&2 \\ \hline
18&6&2&18&8&2&18&9&4&18&18&4 \\ \hline
19&19&2&20&6&2&20&10&2&20&20&3 \\ \hline
21&7&3&21&21&2&22&11&3&22&22&2 \\ \hline
23&23&2&24&6&2&24&8&3&24&12&2 \\ \hline
24&24&3&26&13&2&26&26&2&28&7&3 \\ \hline
28&14&2&28&28&2&30&6&2&30&10&2 \\ \hline
30&15&3&30&30&3&32&8&2&32&16&2 \\ \hline
32&32&2&34&17&2&36&9&2&36&12&2 \\ \hline
36&18&2&36&36&2&38&19&2&40&40&2 \\ \hline
42&7&3&42&14&2&42&21&2&42&42&2 \\ \hline
44&11&2&48&8&2&48&16&2&48&48&2 \\ \hline
54&54&2&60&30&2&60&60&2&66&11&2\\ \hline
70&7&2&84&7&2&&&&&& \\ \hline
\end{tabular}
\caption{Norm List}
\end{table}
\subsection{Discriminant Method $(r \geq 3)$}

This is a refinement of the method used in \cite{MM}, 
and we apply it when  $r \geq 3$. 

If $\Delta$ is the 
discriminant of the power basis $1, \gamma, \gamma^2, \ldots , \gamma^{r-1}$ 
over
$L$ and $\delta_{\IQ(\gamma) \mid L}$, the relative discriminant, then 
$$| N_{L \mid \IQ}(\Delta) | \geq | N_{L \mid \IQ}(\delta_{\IQ(\gamma) \mid L}) |.$$
Choose embeddings $\tau_1, \tau_2, \ldots , \tau_{\mu}$ of $\IQ(\gamma)$ 
into $\IC$ such that $\tau_i|_L = \sigma_i$. Then $N_{L \mid \IQ}(\Delta)
= \prod_{i=1}^{\mu} \sigma_i(\Delta)$ and $\sigma_i(\Delta)$ is the 
discriminant of the power basis $1, \tau_i(\gamma), \tau_i(\gamma^2), \ldots
, \tau_i(\gamma^{r-1})$ of $\tau_i(\IQ(\gamma))$ over $L$. As in 
\cite{MM}, we use Schur's bound \cite{S} which gives that, 
if $-1 \leq x_1 < x_2 < \cdots < x_r \leq 1$ with $r \geq 3$ then
\begin{equation}\label{eqn31}
 \prod_{1 \leq i < j \leq r} (x_i - x_j)^2 \leq M_r = \frac{2^2\,3^3\, \ldots r^r\, 2^2 \, 3^3 \, \ldots (r-2)^{r-2}}{3^3\, 5^5 \, \dots (2r-3)^{2r-3}}.
\end{equation}
Thus, for $i \geq 2$ we have 
 \begin{equation}\label{discr1}
 | \sigma_i(\Delta) | \leq \left( \frac{\sigma_i(\beta_1 \beta_2)}{8} \right)^{r(r-1)} \, M_r.
 \end{equation}
In the case where $i=1$, $\gamma$ has $r-2$ real conjugates over $L$ denoted
by $x_3, x_4, \ldots , x_r$ which, by (\ref{eqn10}), all lie in the interval $(-\beta_1 \beta_2
/4, 0)$. Thus 
\begin{equation}\label{discr2}
| \Delta | \leq | \gamma - \bar{\gamma} |^2 \left(\prod_{i=3}^r (\gamma - x_i)^2(\bar{\gamma}-x_i)^2\right) \left( \frac{\beta_1 \beta_2}{8}\right)^{(r-2)(r-3)} M_{r-2}.
\end{equation}

For $\gamma$ on the contour $\Omega_{p,q}$ we have (see (23))
$$ | \gamma - x_i | = | \bar{\gamma} - x_i | < 4(1+ \cos \pi/p \cos \pi/q)^2.$$ 

We thus define 
$$K_1(p,q,r) = 4 M_{r-2} [4(1+c(p,q))^2]^{4(r-2)} (2 s(p,q)^2)^{(r-2)(r-3)} 
          {\rm Max}_{0 \leq t \leq 1}| \Im(\Omega_{p,q}(t) |^2$$
which can be determined using (\ref{Omega}).
 
 From (\ref{discr1}) and (\ref{discr2}) we obtain an upper bound for $| N_{L \mid \IQ}(\delta_{\IQ(\gamma) \mid L})|$. This is bounded below by 1 but since
$|N_{L \mid \IQ}(\delta_{\IQ(\gamma) \mid L})| = |\Delta_{\IQ(\gamma)}|/\Delta_L^r,$
this lower bound may be improved. 
Since $\IQ(\gamma)$ is a field of degree $r \mu$ with exactly one 
complex place,  
for $n \geq 2$, let $D_n$ denote the minimum absolute value of the discriminant of 
any field of degree $n$ over $\IQ$ with exactly one complex place. For small
values of $n$ the number $D_n$  has been widely investigated (\cite{CDO,Di,DO})
and lower bounds for $D_n$ for all $n$ can be
computed (\cite{Mull,Od,Rodgers,Stark}). In \cite{OdlyzkoU}, the bound 
is given in the form $D_n > A^{n-2} B^2 \exp(-E)$ for varying values of 
$A,B$ and $E$. Choosing, by experimentation, suitable values from this table 
we obtain the bounds shown in Table 3.

\begin{table}[h]
\begin{center}
\begin{tabular}{clcl}
Degree $n$ & Bound & Degree $n$ & Bound  \\
2 & 3 &
3 & 27 \\
4 & 275 &
5 & 4511 \\
6 & 92779 &
7 & 2306599 \\
8 & 68856875* &
9 & $0.11063894 \times 10^{10} $ \\
10 & $0.31503776 \times 10^{11}$ &
11 & $0.90315026 \times 10^{12}$ \\
12 & $0.25891511 \times 10^{14}$ &
13 & $0.74225785 \times 10^{15}$ \\
14 & $0.21279048 \times 10^{17}$&
15 & $0.61002775 \times 10^{18}$ \\
16 & $0.17488275 \times 10^{20}$ &
17 & $0.50135388 \times 10^{21}$ \\
18 & $0.14372813 \times 10^{23} $ &
19 & $0.41203981 \times 10^{24}$ \\
20 & $0.11812357 \times 10^{26}$ 
\end{tabular}
\caption{Discriminant Bounds}
\end{center}
* {\it The exact bound in degree 8 is only known for imprimitive fields \cite{CDO}. This suffices here as the only case not covered here is $p=q=6$ where, by the Norm
List, the degree does not exceed 5.}
\end{table}

For any integer $M \geq 2$, let $D(M) = M^{\phi(M)/2}/(\prod_{\pi} \pi^{\phi(M)/(2 \pi - 2)})$ where the product is over all primes  which divide $M$.
Then
\begin{equation}\label{eqn34}
\Delta_{\IQ(\cos 2 \pi/M)} = \left\{ \begin{array}{ll}
               D(M) & {\rm if~}M \neq m^{\alpha}, 2 m^{\alpha}, m~{\rm a~prime} \\
               D(M)/\sqrt{m} & {\rm if~}M=m^{\alpha}, 2 m^{\alpha},m~{\rm an~odd~prime} \\
               D(M)/2  & {\rm if~}M=2^{\alpha}, \alpha \geq 2.
                \end{array}
         \right.
\end{equation}
If $(p,q) > 2$, 
$L =  \IQ(\cos 2 \pi/M)$ where $M$ is the least common multiple of $p$ and $q$.
If $(p,q) \mid 2$, then 
$\Delta_L = \Delta_{\IQ(\cos 2 \pi/p)}^{\phi(q)/2} \, \Delta_{\IQ(\cos 2 \pi/q)}^{\phi(p)/2}$.

 Thus, from (\ref{discr1}) and (\ref{discr2}), for all $(p,q,r)$ with $r \geq 3$, the following
 inequality must hold
 \begin{equation}{\label{discr}}
 K_1(p,q,r) (2 s(p,q)^2)^{-r(r-1)} \left( \delta_{p,q}/8 \right)^{\mu r(r-1)}
 M_r^{\mu - 1} \geq {\rm Max} \{ 1, D_{r \mu}/\Delta_L^r \}.
 \end{equation}
 Extracting the cases with $r \geq 3$ from the Norm List, and applying this 
 inequality first with a lower bound of 1, results in triples $(p,q,r)$
 where the total degree $r \mu$ is no greater than 20. On these we can apply 
 (\ref{discr}) with values of $D_n$ in Table 3. The result is the so-called {\it 
 Discriminant List} given in Table 4.
\begin{table}[h]
\centering
\begin{tabular}{||c|c|c||c|c|c||c|c|c||}  \hline
p&q&r&p&q&r&p&q&r \\ \hline
6&6&3,4,5&7&7&3,4,5&8&6&3,4 \\ \hline
8&8&3,4,5&9&9&3,4&10&6&3 \\ \hline
10&10&3,4&11&11&3&12&6&3 \\ \hline
12&12&3,4&14&7&3,4&14&14&3 \\ \hline
16&8&3&16&16&3&18&9&3,4 \\ \hline
18&18&3,4&20&20&3&24&8&3 \\ \hline
24&24&3&30&15&3&30&30&3 \\ \hline
\end{tabular}
\caption{Discriminant List}
\end{table}
\subsection{Balancing Method }

Once again this is a refinement of an argument used in \cite{MM}
and here we extend the argument from the case $r=2$ to all $r$. Note 
that the upper bound for $| N(\gamma)|$ used at (\ref{Norm1}) is attained when the
real conjugates of $\gamma$ cluster at one end of the relevant interval,
in which case, the discriminant of the basis using $\gamma$ will be small. 
This argument aims to balance these by incorporating both the norm
amd the discriminant.

Let the minimum polynomial of $\gamma$ over $L$ have roots $x_1 (=\gamma)$, 
$x_2(=\bar{\gamma})$, $x_3 \ldots ,x_r$. 
Recall that,  for each
$\tau_i : \IQ(\gamma) \rightarrow \IR$ such that $\tau_i|_L = \sigma_i$ we have
$\tau_i(\gamma) \in (-\sigma_i(\beta_1 \beta_2/4), 0)$. For $i = 2, \ldots, 
\mu$, let 
$$\tau_i(x_j) = t_i^{(j)}(-\sigma_i(\beta_1 \beta_2 / 4)), \, j= 1,2 \ldots , r$$
so that $0 < t_i^{(j)} < 1$. 
$$N_{L \mid \IQ}({\rm discr}\{1,\gamma, \gamma^2, \ldots, \gamma^{r-1}\}) = 
(\gamma - \bar{\gamma})^2 \prod_{i=3}^r | \gamma - x_i |^4 \prod_{3 \leq j<k\leq r}(x_j-x_k)^2$$
$$~~~~~~~~~~~~~~~~~~ \prod_{i=2}^{\mu}\left( \sigma_i(\frac{\beta_1 \beta_2}{4})^{r(r-1)} \prod_{1 \leq j < k \leq r}(t_i^{(j)} - t_i^{(k)})^2 \right) $$
where here, and later, all empty products have the value 1. Define
$$R_{p,q} = \prod_{i=2}^{\mu} \left(\sigma_i(\frac{\beta_1 \beta_2}{4})^2 \right) = \left(\frac{\delta_{p,q}}{4}\right)^{2 \mu} /(4 s(p,q)^2)^2.$$
Thus 
\begin{equation}{\label{Bal1}}
\prod_{i=2}^{\mu} \prod_{1 \leq j < k \leq r}| t_i^{(j)} - t_i^{(k)} | \geq
\frac{{\rm Max} \{1, (D_{r \mu}/\Delta_L^r)^{1/2}\}}{|\gamma - \bar{\gamma}| \prod_{i=3}^r |\gamma - x_i|^2 \prod_{3\leq j<k\leq r}|x_j-x_k| R_{p,q}^{r(r-1)/4}}.
\end{equation}
On the other hand
$$N_{L \mid \IQ}(N_{\IQ(\gamma) \mid L}(\gamma))= |\gamma|^2 \prod_{i=3}^r x_i
\prod_{i=2}^{\mu}\left(-\sigma_i(\beta_1 \beta_2/4)^r \prod_{j=1}^r t_i^{(j)} \right)$$
so that
\begin{equation}{\label{Bal2}}
\prod_{i=2}^{\mu} \prod_{j=1}^r |t_i^{(j)}| \geq \frac{1}{|\gamma|^2 \prod_{j=3}^r|x_j| R_{p,q}^{r/2}}.
\end{equation}
Let us define $t_i^{(0)} = 0$ for $i = 2, \ldots , \mu$ so that the product of 
(\ref{Bal1}) and (\ref{Bal2}) yields
$$\prod_{i=2}^{\mu} \prod_{0 \leq j<k \leq r}|t_i^{(j)} - t_i^{(k)}| \geq 
\frac{{\rm Max} \{1, (D_{r \mu}/\Delta_L^r)^{1/2}\}}{|\gamma - \bar{\gamma}| |\gamma|^2 \prod_{j=3}^r|\gamma- x_j|^2 \prod_{j=3}^r|x_j| \prod_{3 \leq j < k \leq r}|x_j-x_k| R_{p,q}^{r(r+1)/4}}.$$
Note that $\ds{\prod_{0 \leq j < k \leq r}|t_i^{(j)}-t_i^{(k)}| \leq \left( M_{r+1}/2^{r(r+1)} \right)^{1/2}}$. In the same way, for $r>2$,
$$\prod_{j=3}^r|x_j| \prod_{3 \leq j < k \leq r}|x_j - x_k| \leq \left(M_{r-1} (\frac{\beta_1 \beta_2}{8})^{(r-1)(r-2)}\right)^{1/2}.$$
Also with $\gamma \in \Omega_{p,q}$, $|\gamma-\bar{\gamma}| |\gamma|^2 \prod_{j=3}^r|\gamma - x_j|^2$
will be a maximum when all $x_j$ lie at the left hand extremity of the interval 
$(-\beta_1 \beta_2/4, 0)$. So define
$$K_2(p,q,r) = M_{r-1}^{1/2} (2 s(p,q)^2)^{(r-1)(r-2)/2} \times $$
$$~~~~~~~~~~~~~~~~~~~~~~{\rm Max}_{0 \leq t \leq 1} |2 \Im(\Omega_{p,q}(t)) \Omega_{p,q}(t)^2 (\Omega_{p,q}(t) + 4 s(p,q)^2)^{2(r-2)}|$$
when $r>2$ and $K_2(p,q,2) = {\rm Max}_{0 \leq t \leq 1}| 2 \Im(\Omega_{p,q}(t)) \Omega_{p,q}(t)^2 |$. Thus all our triples $(p,q,r)$ must satisfy
\begin{equation}
K_2(p,q,r) R_{p,q}^{r(r+1)/4} \left(\frac{M_{r+1}}{2^{r(r+1)}}\right)^{(\mu-1)/2} \geq {\rm Max} \{ 1,\left( \frac{D_{r \mu}}{\Delta_L^r} \right)^{1/2}\}.
\end{equation}
Apply this to the Discriminant List for $r \geq 3$ and to the pairs $(p,q)$ 
appearing in the Norm List for $r=2$. In the latter case, if we apply the lower bound of 1 initially, the remaining fields all have total degree not
exceeding 20 and we can then utilise Table 3. The end result is shown 
in Table 5, and termed the {\it Aspiring List}.
\begin{table}
\begin{center}
\begin{tabular}{||c|c|c||c|c|c|||c|c|c||}
\hline
p &q& r & p & q & r& p & q & r\\ 
\hline
6&6&2,3,4,5 &7 &6 & 2 & 7 & 7 & 2,3,4 \\ \hline
8 & 6 & 2,3 & 8 & 8 & 2,3,4,5 & 9 & 6 & 2 \\ \hline
9 & 9 & 2,3  & 10 & 6 & 2,3 & 10  & 10 & 2,3 \\ \hline 
11 & 11 &  2 & 12 & 6 & 2,3 & 12 & 12 & 2,3,4 \\ \hline
13 & 13 &  2  & 14 & 7 & 2,3 & 14 & 14 & 2 \\ \hline
15 & 15 &  2 & 16 & 8 & 2 & 16 & 16 & 2  \\ \hline
18 & 6 &  2  &18 & 9 & 2,3 & 18 & 18 & 2,3 \\ \hline
20 & 10 &  2 & 20 & 20 & 2 &22 & 11 & 2  \\ \hline
24 & 8 &  2 & 24 & 12 & 2 & 24 & 24 & 2  \\ \hline 
28 & 7 &  2 & 30 & 10 & 2 & 30 & 15 & 2  \\ \hline
30 & 30 & 2,3    & 36 & 36 & 2 & 42 & 7 & 2 \\ \hline
42 & 42 &  2  & & & & & & \\ \hline
\end{tabular}
\caption{Aspiring List}
\end{center}
\end{table}

\section{Using the field $L = \IQ(\cos 2 \pi/p, \cos 2 \pi/q)$}

From what we have found so far, the Aspiring List, Table 5, has the following property: 

\medskip

\noindent{\em 
If $\gamma \in \IC \setminus \IR$ is a parameter corresponding to an
arithmetic Kleinian group \\ $\G = \langle  f,g \rangle $ with $f$ of order $p$ and $g$ of order $q$ and  
$[\IQ(\gamma) : L] = r$, then $(p,q,r)$ must appear on the Aspiring List.}

\medskip

Furthermore, $\gamma$ will be an algebraic integer which satisfies an irreducible polynomial 
\begin{equation}\label{eqn35}
x^r + c_{r-1}x^{r-1} + \cdots + c_0 = 0 \quad c_j \in R_L.
\end{equation}
The coefficients $c_j$ are symmetric polynomials in $\gamma, \bar{\gamma}$
and their real conjugates over $L$.  Also the images $\sigma_i(c_j)$
for the real embeddings
$\sigma_i : L \rightarrow \IR$ are symmetric polynomials in the real conjugates
$\tau(\gamma)$ where $\tau : \IQ(\gamma) \rightarrow \IR$ with 
$\tau |_L = \sigma_i$, $i \geq 2$. 

Thus the bounds  on $| \gamma |^2$ and 
$\Re(\gamma)$ obtained from (\ref{Omega}) \S 3 using the freeness criteria and the bounds on the real conjugates $\tau(\gamma)$ in \S 2 using the ramification criteria will 
place bounds on the algebraic integers $c_j$ and $\sigma_i(c_j)$. For each
$(p,q)$ we can readily obtain an integral basis for $L$ over $\IQ$. The bounds on $\gamma$ and its conjugates then translate into bounds on the rational integer
coefficients when each $c_j$ is expressed in terms of this integral basis. 
Once a finite number of possibilities for each coefficient $c_j$ individually
is obtained, the roots of each of the resulting finite number of polynomials 
at (\ref{eqn35}) so obtained, and their conjugates, can be further examined to 
see if their roots satisfy the required bounds. We explain the basic methods used to carry out this 
computational process in this section.  This basic method is carried 
out as a first step by a simple Maple program on the triples in the Aspiring List.

\bigskip

These remarks above actually apply to any algebraic integer $\delta$ in $\IQ(\gamma)$
such that $\IQ(\delta) = \IQ(\gamma)$ and for which one can obtain bounds on 
$\delta$ and its conjugates. In particular, if $v$ is a unit in $L$,
we can take $\delta = \gamma/v$
and  suitable choices of $v$ lead to improved bounds
on $\delta$. 

\bigskip

For the basic method which we now describe, we assume first that $\mu \geq 3$, the
cases where $\mu \leq 2$ being considerably easier to handle.
For all $(p,q,r)$ on the Aspiring List,  $L$ has an integral 
basis of the form $\{ 1, u, u^2, \ldots , u^{\mu - 1} \}$ where 
$u = 2 \cos 2 \pi/M$ for some integer $M$. 
 Let $\sigma_1 = {\rm Id}, \sigma_2, \ldots , \sigma_{\mu}$ denote
the Galois automorphisms of $L$ over $\IQ$ with $\sigma_i (2 \cos 2 \pi/M) 
= 2 \cos 2 \pi y_i/M$ where $1 \leq y_i < M/2$ and $(y_i,M)=1$.

Let $\delta$ be an algebraic integer as described above which 
satisfies (\ref{eqn35}).
Let 
\[ c_j = m_0 + m_1 u + m_2 u^2 + \cdots + m_{\mu-1} u^{\mu-1}\]
 where $m_k \in \IZ$.
Let $A$ be the $\mu \times \mu$ matrix $[\sigma_i(u^{j-1})]$, $1 \leq i,j \leq \mu$. Then 
\begin{equation}\label{eqn36} A (\tilde{m}) = \tilde{c_j} \end{equation}
where $\tilde{m} = (m_0, m_1, \ldots , m_{\mu-1})^t$ and $\tilde{c_j} = 
(c_j, \sigma_2(c_j), \ldots , \sigma_{\mu}(c_j))^t$. Thus 
\begin{equation}\label{eqn37}
\tilde{m} = A^{-1} \tilde{c_j}
\end{equation}
where we can numerically determine the entries of $A$ and $A^{-1}$. 
The bounds on $| \gamma |^2, \Re(\gamma)$ obtained by maximising them on 
$\Omega_{p,q}$ using (\ref{Omega}) and the bounds on the real 
conjugates at (\ref{eqn10})
give bounds on $\delta = \gamma/v, v \in R_L^*$ and its conjugates and hence
on each entry of the matrix $\tilde{c_j}$.
 Thus there exist
$\mu \times 1$ matrices $I_j$ and $S_j$ such that $I_j \leq \tilde{c_j} \leq S_j$
with the obvious notation. In the cases where $p$ is not a prime power, 
$4 \sin^2 \pi/p = -\beta_1$ is a unit and in these cases it is expedient
 to take 
$\delta = \gamma/(-\beta_1)$ or $= \gamma/(-\beta_2)$ if $q$ is 
also not a prime power. 

\begin{example} 
$(p,q,r) = (42,42,2)$. In this case with $\delta = \gamma/(-\beta_1)$,
$c_0 = | \gamma |^2/(16 \sin^4 \pi/42)$ and $0 < \sigma_i(c_0) < 
\sigma_i(\beta_1 \beta_2 / 4 \beta_1)^2 = \sigma_i(\sin^2 \pi/42)^2$.
Thus $I_0 < \tilde{c_0} < S_0$ with $I_0 = 0$, and the $i$-th entry $s_i$ of $S_0$ is 
$\sigma_i(\sin^2 \pi/42)^2$ for $i=2,3, \ldots ,6$ and $s_1 = 16(\cos^2 \pi/42/
\sin^2 \pi/42)^2$.
\end{example}

\noindent{\bf Remark.} From this example, a common feature of many examples will be noted - that
all entries of $S_0$ except the first are small. This is a consequence of our 
 choice of $v$ and we will explain below how to exploit this.
 
 \bigskip

Let us return to the general case as at (\ref{eqn37}). We can obtain upper and
 lower estimates on $\tilde{m}$ as follows:  Write $A^{-1} = A_+^{-1} + A_-^{-1}$ where 
$A_+^{-1}, A_-^{-1}$ are $\mu \times \mu$ matrices with all entries in 
$A_+^{-1}$ being $\geq 0$ and those in $A_-^{-1}$ being $\leq 0$. 
 We thus obtain
\begin{equation}\label{eqn38}
A_+^{-1} I_j + A_-^{-1} S_j \leq \tilde{m} \leq A_+^{-1} S_j + A_-^{-1} I_j.
\end{equation}
This then gives a finite number of possibilities for $\tilde{m}$. We 
refer to this as a search space and from these inequalities,  its size can be 
readily measured. In general, the search space described by (\ref{eqn38}) can be extremely
large.
In 
Example 5.1 above,  for example, 
 it is of the order of $1.5 \times 10^{25}$. In such cases,  we extend this 
technique to exploit the fact that, in many cases, all the entries of $I_j$ and 
$S_j$ except the first are small.

From (\ref{eqn38}) determine the possible values of $m_0$, the first entry of $\tilde{m}$
and the constant term in the expression of $c_j$ in terms of the integral basis
$1, u, u^2, \ldots, u^{\mu -1}$. For each $m_0$ we have 
\begin{equation}\label{eqn39}
m_1 u + m_2 u^2 + \cdots + m_{\mu - 1} u^{\mu - 1} = c_j - m_0
\end{equation}
and the corresponding $\mu - 1$ equations under the embeddings $\sigma_i, 
i = 2, \ldots, \mu$. Now if $B$ denotes the $\mu - 1 \times \mu - 1 $ matrix 
obtained from $A$ by deleting the first row and first column and if 
$\tilde{m}', \tilde{c_j}'$ denote the $\mu - 1 \times 1$ matrices obtained 
by removing the first entries of $\tilde{m}, \tilde{c_j}$, 
we can write the $\mu-1$ equations obtained from (\ref{eqn39}) for the
 embeddings $\sigma_2, \ldots, \sigma_{\mu}$, in the form
$$
B \tilde{m}' = \tilde{c_j}' - m_0 \tilde{1}
$$
where $\tilde{1}$ is the $\mu - 1 \times 1$ matrix all of whose entries are 1.
This then yields $\tilde{m}' = B^{-1} \tilde{c_j}' - m_0 B^{-1} \tilde{1}$.
For each $m_0$ the term $m_0 B^{-1} \tilde{1}$ is fixed. By splitting 
$B^{-1}$ into its positive and negative parts as we did for $A^{-1}$ and using 
the truncated limits ${I_j}',{S_j}'$ for $\tilde{c_j}$  we obtain bounds on $\tilde{m}'$ given by 
\begin{equation}\label{eqn40}
B_+^{-1} {I_j}' + B_-^{-1} {S_j}' - m_0 B^{-1} \tilde{1} \leq \tilde{m}' \leq B_+^{-1} {S_j}' + B_-^{-1} {I_j}' - m_0 B^{-1} \tilde{1}.
\end{equation}
If the entries of ${I_j}', {S_j}'$ are small, this yields a small search space for 
$\tilde{m'}$ whose size is essentially independent of $m_0$.
In Example 5.1, for example,  there are 6166 possibilities for $m_0$
and 576 for $\tilde{m}'$ so that the search space is now of the 
order of $3.5 \times 10^6$,  a significant reduction.

 For each resulting $\tilde{m}$ we check the validity
of $I_j \leq A \tilde{m} \leq S_j$ and list the resulting $\tilde{m}$ and hence
candidate $c_j$. Again in Example 5.1,  there are three  such integer vectors 
$\tilde{m}$ and hence only three candidates for $c_0$.

\medskip

In the cases where $\mu=2$, we dispense with the use of the matrix $A$ 
(and hence B). For in that case, all integers in $R_L$ have the form 
$(a+b \sqrt{d})/2$ where $a,b \in \IZ$ with $a \equiv b ({\rm mod}~2)$ 
and $a \equiv b \equiv 0 ({\rm mod}~2)$ if $d \not\equiv 1 ({\rm mod}~4)$.
Thus if $c_j = (a_j + b_j \sqrt{d})/2$, the upper and lower bounds on 
$c_j$ and $\sigma(c_j)$ respectively for $\sigma$ the non-identity 
embedding, can be expressed as
$$\ell_1 < \frac{a_j+b_j \sqrt{d}}{2} < u_1, \quad \ell_2 < \frac{a_j - b_j \sqrt{d}}{2} < u_2.$$
Thus $a_j$ must be an integer between $\ell_1 + \ell_2$ and $u_1 + u_2$ and, for 
each such $a_j$, $b_j$ lies between $(2 \ell_1 - a_j)/\sqrt{d}$ and 
$(2 u_1 - a_j)/\sqrt{d}$. Provided the bounds are reasonable, it is a simple 
matter to find all the integers satsifying these inequalities.

These methods described above for enumerating and listing candidate values of 
the coefficients $c_j$ in either the cases where $\mu \geq 3$ or $\mu = 2$ 
will be referred to as the {\em Basic Method}.

\medskip

In the (many) cases where $p=q$, we noted in \S 2 that any non-elementary 
Kleinian group generated by $f,g$ where $o(f)=o(g)=p$ is a subgroup of 
index at most 2 in a non-elementary Kleinian group generated by $f$
and an element $h$ of order 2. Thus, in these cases, by Theorem 2.3, instead
of trying to determine $\gamma = \gamma(f,g)$, we can search for possible
values of $\gamma_1 = \gamma(h,f)$. For a real embedding $\tau : \IQ(\gamma) 
\rightarrow \IR$ with $\tau_L = \sigma$ we have, by (\ref{eqn13}), $-\sigma(4 
\sin^2 \pi/p) < \tau(\gamma_1) < 0$. Also from \S 3, $| \gamma_1 | \leq 4$.
Furthermore, by (\ref{eqn11}), $\gamma_2 = \beta_1 - \gamma_1$ also corresponds to a 
group generated by an element of order 2 and an element of order $p$. Thus 
we can assume that the $\gamma_1$-space is symmetric about $\Re(\gamma_1) 
= \beta_1/2$ and so
\begin{equation}
 - 4 < \Re(\gamma_1) < -2 \sin^2 \pi/p.
\end{equation}
We can thus apply the same strategy as in the {\em Basic Method} to determine
the coefficients $c_j$ of the polynomial satisfied by $\gamma_1$ or $\delta_1
= \gamma_1/v$ for a suitable unit $v \in R_L$. We refer to this also as a
{\em Basic Method}.

Applying  the {\em Basic Methods}  to 
 triples on the Aspiring List yields candidate values for the coefficients
$c_j$ of the polynomials $p(x)$ satisfied by some $\delta$ where 
$\IQ(\delta) = \IQ(\gamma)$.
In some cases the bounds are tight enough that there are no candidate
values for one of the coefficients. We list these below in Table 6. 
In this Table and subsequently, we will use the notation $\gamma(p,q)$
for $\gamma(f,g)$ where $o(f)=p, o(g)=q$ and also $\gamma(2,p)$ for 
$\gamma(h,f)$ where $o(h)=2$.
Generally, the search spaces are small in the cases of coefficients 
$c_0$ and $c_{r-1}$ as they are, up to sign, the product and sum of the roots. 
Thus degree 2, considered in \S 6 below, is reasonably straightforward.
For the other coefficients, additional methods may be required to reduce 
the size of the search space to manageable proportions.
These will be discussed in \S 8 to 9 below.
\begin{table}[h]
\begin{center}
\begin{tabular}{|c|c|c|}
\hline
Triple & $\delta$ & Outcome \\ \hline
(28,7,2) & $\gamma(28/7)/4 \sin^2 \pi/28$ & No values of $c_0$ \\ \hline
(22,11,2) & $\gamma(22,11)$ & No values of $c_0$ \\ \hline
(16,8,2) & $\gamma(16,8)/(1+2 \cos 6 \pi/16)$ & No values of $c_1$. \\ \hline
\end{tabular}
\caption{}
\end{center}
\end{table}

\section{Degree 2}
         Here we consider the cases where $r = 2$ so that $\delta$ satisfies
$p(x) = x^2 + c_1 x + c_0$. From the Basic Methods we have obtained 
candidate values for $c_1$ and $c_0$. The polynomial $p(x)$ will define
 a field with one complex place if and only if $c_1^2 - 4 c_0 <0$ and 
 $\sigma_i(c_1^2-4c_0)>0$  for $i = 2,3, \ldots, \mu$. Furthermore, for 
 $i \geq 2$ both roots of $p^{\sigma_i}(x)=x^2 + \sigma_i(c_1)x + \sigma_i(c_0)
 = 0$ must lie in an interval $(-\ell_i, 0)$ where $\ell_i>0$ is the bound
 obtained using (\ref{eqn10}) or (\ref{eqn13}) for the particular
 choice of $\delta$. By the Basic Method, $0 \leq \sigma_i(c_1) < 2 \ell_i$
 and $0 < \sigma_i(c_0) < \ell_i^2$ for $i \geq 2$. Thus the condition on 
 the location of these real roots is equivalent to requiring that 
 $\ell_i^2 - \sigma_i(c_1) \ell_i + \sigma_i(c_0) > 0$. 
 Thus all 
these conditions can be checked directly on the candidate coefficients 
$c_1, c_0$. This will be referred to as {\em polynomial reduction}.
\begin{examples}

\noindent(1.) $(p,q,r) = (42,42,2)$. As in Example 5.1, take $\delta = 
\gamma(42,42)/4 \sin^2 \pi/42$. The Basic Method throws up two 
candidates for $c_1$ and three for $c_0$. None of the 6 resulting 
polynomials satisfy all the inequalities above and so there are no 
arithmetic Kleinian groups corresponding to the triple $(42,42,2)$.

\noindent(2.) $(p,q,r) = (24,24,2)$. Taking $\delta = \gamma(2,24)$
the Basic Method yields 74 candidates for $c_0$ and 20 for $c_1$. Then
polynomial reduction reduces this to two polynomials.

\noindent(3.) $(p,q,r)= (12,6,2)$. With $\delta = \gamma(12,6)$ we obtain 
45 candidates for $c_0$ and 19 for $c_1$ and polynomial reduction reduces 
this to a total of 45 polynomials.
\end{examples}
The remaining polynomials can then be computationally solved and the complex 
roots checked to see if they give rise to values of $\gamma(p,q)$ 
which lie inside the contour $\Omega_{p,q}$.
All the polynomials which are left at this stage correspond to a $\gamma$
which satisfies conditons 1,2 and 3 of Theorems 2.2 or 2.3. If $p=q$
and the deduction is carried out using $\gamma(2,p)$, then the 
resulting $\gamma(p,p)= \gamma(2,p)(\gamma(2,p) + 4 \sin^2 \pi/p)$ 
corresponds to a subgroup
of an arithmetic Kleinian group. It can turn out that the resulting
$\gamma(p,p)$ is real, which cases, as noted in \S 2, are completely
understood.
\begin{example} $(p,q,r)=(24,24,2)$. The two polynomials (see above)
both yield that $\gamma(24,24)$ is real and there are no such arithmetic
Kleinian groups.(See \S 2.3). 
\end{example}

More generally, we still need to check condition 4 of Theorem 2.2 for 
$\gamma(p,q)$ using Lemma 2.4. If $p=q$ and $\gamma(p,p)$ is 
obtained by first determining  $\gamma(2,p)$, then this condition
is automatically satisfied as noted at the end of \S 2.3. Thus this 
is most frequently applied in the cases where $ p \neq q$.
\begin{example} $(p,q,r)=(12,6,2)$. Here the field $L= \IQ(\sqrt{3})$
and we have 45 candidate polynomials from above. Using (\ref{eqn15})
we replace the variable $x$ by $y(y-\sqrt{3})/(3(2+\sqrt{3}))$ and find
that just one of the resulting quartic polynomials in $y$ factorise
in $\IQ(\sqrt{3})$. Thus there is one value of $\gamma(12,6)$ 
which gives rise to a subgroup of an arithmetic Kleinian group in this case.
\end{example}
Using this {\em factorisation method} any remaining polynomials will
give values of $\gamma$ which correspond to subgroups of arithmetic 
groups. The results are shown in Table 7. These parameters must then be subjected to geometric methods to ascertain
if they have finite covolume and so satisfy the final conditions 
of Theorems 2.2 or 2.3. These geometric methods will be described in \S 10.

The notation used in Table 7 is as follows: the second column gives the generating
element $\delta$ to which we apply the Basic Method. The next two columns give
the number of resulting possible values of the coefficients $c_0$ and $c_1$.
The column headed ``PR'', refers to the 
number of polynomials remaining after polynomial reduction, that headed 
``B'' gives the number that are non-real and lie inside the contour
$\Omega_{p,q}$ and the ``F'' column those left after the factorisation 
criteria has been applied. Thus the non-zero entries in the final column
are those which need to be further considered by geometric methods.(The * 
in the $(42,7,2)$ row indicates that the values of $c_1$ were calculated
and from the small number of resulting values we obtained improved bounds
on $c_0$ by using the inequalities implied by the method of polynomial
 reduction. The - in the row of $(14,7,2)$ indicates that we omitted this step.)
\begin{table}
\begin{center}
\begin{tabular}{|c|c|c|c|c|c|c|} \hline
Triple & $\delta$ & $c_0$ & $c_1$ & PR & B & F \\ \hline
(42,42,2) & $\gamma(42,42)/4 \sin^2 \pi/42$ & 3 & 2 & 0 & 0 & 0 \\ \hline
(42,7,2) & $\gamma(42,7)/4 \sin^2 \pi/42 \times 4 \sin^2 \pi/21$ & * & 5 & 0 & 0 & 0\\ \hline
(36,36,2) & $\gamma(2,36)$ & 16 & 10 & 0 & 0 & 0 \\ \hline
(30,30,2) & $\gamma(2,30)$ & 249 & 44 & 10 & 1 & 1 \\ \hline
(30,15,2) & $\gamma(30,15)/4 \sin^2 \pi/15$ & 36 & 20 & 0 & 0 & 0 \\ \hline
(30,10,2) & $\gamma(30,10)/4 \sin^2 pi/10$ & 9 & 8 & 0 & 0 & 0 \\ \hline
(24,24,2) & $\gamma(2,24)$ & 72 & 20 & 2 & 0 & 0 \\ \hline
(24,12,2) & $\gamma(24,12)/4 \sin^2 \pi/12$ & 6 & 5 & 0 & 0 & 0 \\ \hline
(24,8,2) & $\gamma(24,8)$ & 12 & 12 & 1 & 0 & 0 \\ \hline
(20,20,2) & $\gamma(20,20)/4 \sin^2 \pi/20$ & 16 & 13 & 0 & 0 & 0 \\ \hline
(20,10,2) & $\gamma(20,10)/4 \sin^2 \pi/20$ & 1 & 4 & 0 & 0 & 0 \\ \hline
(18,18,2) & $\gamma(2,18)$ & 122 & 30 & 16 & 3 & 3 \\ \hline
(18,9,2) & $\gamma(18,9)/4 \sin^2 \pi/18$ & 268 & 62 & 73 & 47 & 2 \\ \hline
(18,6,2) & $\gamma(18,6)/4 \sin^2 \pi/18$ & 6 & 9 & 0 & 0 & 0 \\ \hline
(16,16,2) & $\gamma(2,16)$ & 61 & 19 & 0 & 0 & 0 \\ \hline
(15,15,2) & $\gamma(15,15)/4 \sin^2 \pi/15$ & 4 & 5 & 0 & 0 & 0 \\ \hline
(14,14,2) & $\gamma(14,14)/4 \sin^2 \pi/14$ & 85 & 38 & 10 & 3 & 0 \\ \hline
(14,7,2) & $\gamma(14,7)/4 \sin^2 \pi/14$ & 244 & 65 & 161 & - & 1 \\ \hline
(13,13,2) & $ \gamma(2,13)$ & 11 & 13 & 0 & 0 & 0 \\ \hline
(12,12,2) & $ \gamma(2,12)$ & 64 & 17 & 67 & 18 & 18 \\ \hline
(12,6,2) & $\gamma(12,6)$ & 45 & 19 & 45 & 30 & 1 \\ \hline
(11,11,2) & $\gamma(2,11)$ & 35 & 17 & 0 & 0 & 0 \\ \hline
(10,10,2) & $\gamma(2,10)$ & 48 & 8 & 44 & 15 & 15 \\ \hline
(10,6,2) & $ \gamma(10,6)$ & 34 & 20 & 40 & 24 & 0 \\ \hline
(9,9,2) & $ \gamma(2,9)$ & 72 & 22 & 7 & 2 & 2 \\ \hline
(9,6,2) & $\gamma(9,6)$ & 4 & 7 & 1 & 0 & 0 \\ \hline
(8,8,2)& $ \gamma(2,8)$ & 65 & 17 & 48 & 20 & 20 \\ \hline
(8,6,2) & $\gamma(8,6)$ & 42 & 21 & 50 & 33 & 0 \\ \hline
(7,7,2) & $\gamma(2,7)$ & 199 & 43 & 32 & 8 & 8 \\ \hline
(7,6,2) & $\gamma(7,6)$ & 8 & 12 & 0 & 0 & 0 \\ \hline
(6,6,2) & $\gamma(2,6)$ & 16 & 8 & 78 & 24 & 24 \\ \hline
\end{tabular}
\caption{Degree 2 candidates}
\end{center}
\end{table}

\section{Degree 3}.
Apart from the {\em polynomial reduction} process, this is very similar to
the degree 2 cases as carried out in the preceding section. Let $\delta$ be 
such that $\IQ(\delta) = \IQ(\gamma)$ where $[Q(\gamma) : L ] = 3$ so that
$\delta$ satisfies $p(x) = x^3 + c_2 x^2 + c_1 x + c_0=0$ with $c_i \in R_L$.
Using the Basic Methods we obtain candidate values for $c_0, c_1$ and $c_2$.
In general, there are many more candidates for $c_1$ than for $c_0$ or $c_2$.
We then ascertain that at the non-identity real places $\sigma_i$ of $L$,
the conjugate polynomials $p^{\sigma_i}(x)$ has three real roots in the interval
$(-\ell_i,0)$ where $\ell_i$ is obtained from (10) and (13). This can be checked
without numerically solving the polynomial (which is a time consuming process)
by the following sequence of requirements on combinations of the coefficients:
\begin{itemize}
\item $\sigma_i(c_2)^2 > 3 \sigma_i(c_1)$, which forces the derivative 
$Dp^{\sigma_i}(x)$ to have two real roots;
\item $Dp^{\sigma_i}(-\ell_i) > 0$, which forces these roots, $r_1,r_2$ to lie 
in the interval $(-\ell_i, 0)$;
\item $p^{\sigma_i}(-\ell_i) < 0$, which forces $p^{\sigma_i}(x)$ to have at least one root in the interval $(-\ell_i,0)$;
\item $\sigma_i(-2 c_2 c_1 c_0/3 + 4 c_1^3/27+4 c_2^3 c_0/27-c_2^2 c_1^2/27 + c_0^2)<0$. which forces $p^{\sigma_i}(r_1)p ^{\sigma_i}(r_2)<0$ and so $p^{\sigma_i}(x)$ to have three real roots in the interval $(-\ell_i,0)$.
\end{itemize}
Any cubic remaining after this, can then be solved 
at the identity real place of $L$ to ensure that it has a pair of non-real 
roots and that the real root lies in the interval $(-\ell_1,0)$. 
Following this {\em polynomial reduction} procedure, we check to determine
if the values of $\gamma(p,q)$ lie inside the contour $\Omega_{p,q}$.
Finally, if appropriate, we apply the factorisation condition of 
Theorem 2.2. The results are tabulated in Table 8, 
as in Table 7 so that any non-zero 
numbers in the right hand column correspond to groups which must be checked
by geometric methods to see if they have finite covolume.

\begin{table}
\begin{center}
\begin{tabular}{|c|c|c|c|c|c|c|c|} \hline
Triple & $\delta$ & $c_0$ & $c_1$ & $c_2$ & PR & B & F \\ \hline
(30,30,3) & $\gamma(2,30/4 \sin^2 pi/30)$ & 250 & * & 296 & - & - & 0 \\ \hline
(18,18,3) & $\gamma(2,18)/4 \sin^2 pi/18$ & 8 & 4442 & 180 & 1 & 0 & 0 \\ \hline
(18,9,3) & $\gamma(18,9)/4 \sin^2 \pi/18$ & 11 & 2429 & 137 & 0 & 0  & 0 \\ \hline
(14,7,3) & $\gamma(14,7)/4 \sin^2 \pi/14$ & 25 & 2207 & 148 & 1 & 0 & 0  \\ \hline
(12,12,3) & $ \gamma(2,12)$ & 65 & 218 & 26 & 85 & 19 & 19 \\ \hline
(12,6,3) & $\gamma(12,6)$ & 3 & 138& 30 & 1 & 1 & 1 \\ \hline
(10,10,3) & $\gamma(2,10)$ & 48 & 175 & 29 & 33 & 5 & 5  \\ \hline
(10,6,3) & $ \gamma(10,6)$ & 1 & 103 & 32  & 0 & 0 & 0 \\ \hline
(9,9,3) & $ \gamma(2,9)$ & 219 & 812 & 56 & 0 & 0 & 0 \\ \hline
(8,8,3)& $ \gamma(2,8)$ & 133 &   256&   30& 268 & 29 & 29 \\ \hline
(8,6,3) & $\gamma(8,6)$ & 5 & 129 &  32& 2 & 0 & 0\\ \hline
(7,7,3) & $\gamma(2,7)$ & 1381 & 2449 & 105 & 26 & 1 & 1  \\ \hline
(6,6,3) & $\gamma(2,6)$ & 16 & 24 & 9 & 1496 & 124 & 124 \\ \hline
\end{tabular}
\caption{Degree 3 candidates}
\end{center}
\end{table}
Note: The * in case $(30,30,3)$ indicates that we actually used the 
linked triples method which is explained in the next sections. The outcome 
was that there were no linked triples and so no groups can arise.

\section{Degree $ \geq 4$}

From the {\it Aspiring List}, we see that there are six cases with $ r \geq  4$
all with $p=q$. In general, the Basic Methods enable one to determine the 
candidates for the coefficients $c_0$ and $c_{r-1}$, but give rise to 
unfeasible search spaces in attempting to determine the other coefficients.
So we develop some new methods of obtaining bounds on the coefficients 
by exploiting the relationship at (11) between $\gamma = \gamma(p,p)$
and $\gamma_1 = \gamma(2,p)$. Since $\gamma_2 = \beta_1 - \gamma_1$ is also
a candidate $\gamma(2,p)$ value, equation (11) can be stated as 
\begin{equation}
\gamma = - \gamma_1\, \gamma_2 .
\end{equation}
We use the most ``awkward'' case $(7,7,4)$, which is the case of highest
total degree over $\IQ$ amongst these six, as a template to describe our methods.

Let $\beta_1 = -(2 - 2 \cos 2 \pi/7)$ and $\beta_i = \sigma_i(\beta_1)$
where $\sigma_i$, $i= 2,3$ are the non-trivial automorphisms of $L = 
\IQ(\cos 2 \pi/7)$. Let $B_i = - \beta_i^2/4$ so that, for $\tau : k \Gamma 
\rightarrow \IR$, we have 
\begin{equation}
\beta_i < \tau(\gamma_1), \tau(\gamma_2) < 0 ~{\rm and~} B_i < \tau(\gamma)
< 0
\end{equation}
where $\tau |_L = \sigma_i$. From \S3, we also have bounds on the complex number
$\gamma$ i.e.
\begin{equation}
| \gamma | < 4(2 \cos \pi/7)^2 = G_u ~{\rm and~} R_{\ell} < \Re(\gamma) < G_u
\end{equation}
where $R_{\ell} \approx -5.0914 $ computed using (19). Since $\gamma_1, \gamma_2$
are symmetric about $\Re(\gamma_i) = \beta_1/2$, we can assume that 
$| \gamma_1 | <4$ and $-4 < \Re(\gamma_1) \leq \beta_1/2 $ and $|\gamma_2 | < 4$
and $\beta_1/2 \leq  \Re(\gamma_2) < \beta_1 + 4$ ( see \S 3).

Let $\gamma, \gamma_1, \gamma_2$ satisfy the polynomials
\begin{equation}\left.
\begin{array}{ccl}
p(x) &= &  x^4 + c_3 x^3 + c_2 x^2 + c_1 x + c_0 \\  
p_1(x)& = & x^4 + c_3^{(1)} x^3 + c_2^{(1)} x^2 + c_1^{(1)} x + c_0^{(1)} \\ 
p_2(x)& = & x^4 + c_3^{(2)} x^3 + c_2^{(2)} x^2 + c_1^{(2)} x + c_0^{(2)}
\end{array} \right\}
\end{equation}
respectively. 
As noted above, we can determine $c_0, c_0^{(1)}, c_0^{(2)}, c_3, c_3^{(1)}, 
c_3^{(2)}$ by our Basic Methods. The basic ideas here are then to use 
these determined values to place bounds and restrictions on the remaining
coefficients. Furthermore, since $\gamma_2 = \beta_1 - \gamma_1$, the 
coefficients of $p_2(x)$ are combinations of the coefficients of $p_1(x)$.
All this enables us to determine $c_2^{(1)}$ and $c_1^{(1)}$ from the other 
coefficients.

Using the Basic Methods we determine candidates for $c_0$ and $c_0^{(1)}$ 
(and hence $c_0^{(2)}$). There are 412 and 9769 respectively. From (42)
it follows that $c_0 = c_0^{(1)} c_0^{(2)}$ and we determine
 all such linked triples $(c_0, c_0^{(1)}, c_0^{(2)} )$. (In this 
 $(7,7,4)$ case it is expedient to first narrow down the search by using the 
 fact that the rational integral equation $N_{L \mid \IQ}(c_0) =
 N_{L \mid \IQ}(c_0^{(1)}) N_{L \mid \IQ}( c_0^{(2)})$ must hold.)
 There are 8979 linked triples. Since $\gamma_2 = \beta_1 - \gamma_1$, then
 \begin{equation}
 c_0^{(2)} = \beta_1^4 + c_3^{(1)} \beta_1^3 + c_2^{(1)} \beta_1^2 +
 c_1^{(1)} \beta_1 + c_0^{(1)} .
 \end{equation}
 This implies that $\beta_1 \mid c_0^{(1)}- c_0^{(2)} $, which, if 
 $c_0^{(1)} - c_0^{(2)} = a + bu + c u^2$ where $u = 2 \cos 2 \pi/7$ is 
 equivalent to $a + 2b + 4c  \equiv 0({\rm mod}~7)$. We reduce our set of 
 linked triples to satisfy this divisiblity condition, obtaining 1303 
 such triples.

For each candidate linked triple, we now obtain new bounds on the coefficients 
$c_1, c_1^{(1)}, c_1^{(2)}$ and their conjugates which depend on the values of 
a linked triple as follows:  Let the roots of $p^{\sigma_i}(x), i = 2,3$ 
be $y_1, y_2, y_3, y_4$ so that 
$\sigma_i(c_0) =  y_1 y_2 y_3 y_4$ and $\sigma_i(c_1) = - (y_1 y_2 y_3 y_4) \sum_{j=1}^4(1/y_j)$.
Let the $y_j$ be ordered so that $B_i < y_4 < y_3 < y_2 < y_1 < 0$. So
$\sigma_i(c_0) < (-y_1)(- B_i)^3$ and thus $(-1/y_1) <(- B_i)^3/\sigma_i(c_0)$. Also
$\sigma_i(c_0) < (-y_2)^2(-B_i)^2$ so that $(-1/y_2) < ((-B_i)^2/\sigma_i(c_0))^{1/2}$. Continuing in this vein, we obtain
\begin{equation}\label{eqn57}
\sigma_i(c_1) < \sigma_i(c_0)\left[ \frac{(-B_i)^3}{\sigma_i(c_0)} + \left(\frac{(-B_i)^2}{\sigma_i(c_0)}\right)^{1/2} + \left(\frac{-B_i}{\sigma_i(c_0)}\right)^{1/3} + \left( \frac{1}{\sigma_i(c_0)}\right)^{1/4}\right].
\end{equation}
In the other direction, from the arithmetic/geometric mean inequality, we deduce that
\begin{equation}\label{eqn58}
\sigma_i(c_1) \geq 4 \sigma_i(c_0)^{3/4}.
\end{equation}
In a similar way at the identity embedding,   we obtain an upper bound on $c_1$ as
\begin{equation}\label{eqn59}
c_1 < c_0 \left[ \left(\frac{G_u^2 (-B_1)}{c_0} \right) + \left( \frac{G_u^2}{c_0}\right)^{1/2} \right] - 2 R_{\ell} B_1^2.
\end{equation}
On the other hand,
$$ \frac{c_1}{c_0} = - \left( \frac{1}{x_3} + \frac{1}{x_4} \right) - \left( \frac{1}{\gamma} + \frac{1}{\bar{\gamma}}\right),$$
where the roots of $p(x)$ are $\gamma, \bar{\gamma},x_3,x_4$
The first term here is greater than $-2/B_1$ and the second is greater than
$-2/|\gamma|$.  J{\o}rgensen's Lemma states $|\gamma |+|\beta_1|\geq1$ in a 
discrete non-elementary group, thus $|\gamma | > 2 \cos 2\pi/7 - 1$ so that
\begin{equation}\label{eqn60}
c_1 > c_0 \left( \frac{-2}{B_1} - \frac{2}{2 \cos 2 \pi/7-1} \right).
\end{equation}

In an entirely analogous manner, we can obtain similar bounds for $c_1^{(1)}$ and $c_1^{(2)}$ depending on each $c_0^{(1)}$ and $c_0^{(2)}$ in a linked triple.

Thus for $i=2,3$ and $j = 1,2$ we have
\begin{equation}
 4\sigma_i(c_0^{(j)})^{3/4} < \sigma_i(c_1^{(j)})<\sigma_i(c_0^{(j)})\left[\frac{(-\beta_i)^3}{\sigma_i(c_0^{(j)})}+ \cdots + \left(\frac{1}{\sigma_i(c_0^{(j)})}\right)^{1/4}\right]. 
 \end{equation}
Using the symmetry of $\gamma_1, \gamma_2$, we obtain 
\begin{equation}\label{eqn61}
\left(\frac{-2}{\beta_1}-\frac{\beta_1}{16}\right)c_0^{(1)} < c_1^{(1)} < 
c_0^{(1)} \left[ \left(\frac{16 (-\beta_1)}{c_0^{(1)}}\right) +  \left(
\frac{16}{c_0^{(1)}}\right)^{1/2}\right] + 8 \beta_1^2.
\end{equation}
If $p_2(x)$ has roots $\gamma_2, \bar{\gamma_2}, z_3, z_4$, then
$$c_1^{(2)} = -| \gamma_2 |^2 (z_3 + z_4) - (\gamma_2 + \bar{\gamma_2}) z_3 z_4.$$
Using the AM/GM inequality and the fact that $-1 < \beta_1 < z_3, z_4 <0$ we have
\begin{equation}
c_1^{(2)} \geq | \gamma_2|^2 2(z_3 z_4)^{1/2}-2 \Re(\gamma_2)z_3 z_3 > 2(z_3z_4)(|\gamma_2|^2 - 2 \Re(\gamma_2)) > - \beta_1^2/2.
\end{equation}
Also
\begin{equation}
c_1^{(2)} < c_0^{(2)} \left[\left( \frac{16(-\beta_1)}{c_0^{(1)}}\right) + \left(\frac{16}{c_0^{(1)}}\right)^{1/2}
\right]-\beta_1^3.
\end{equation}

We now further exploit relation (11) to deduce that
\begin{equation}\label{eqn62}
 - \beta_1 c_1 = c_0^{(2)} c_1^{(1)} + c_0^{(1)} c_1^{(2)}.
\end{equation}
This gives upper bounds for $c_1^{(1)}$ and its conjugates which, in many cases, are an improvement on those obtained at (51) and (52) since
\begin{equation}\label{eqn63}
c_1^{(1)} \leq \frac{-\beta_1}{c_0^{(2)}} (\mbox{maximum value of $c_1$}) - \frac{c_0^{(1)}}{c_0^{(2)}}(\mbox{minimum value of $c_1^{(2)}$}).
\end{equation}

In fact, in this $(7,7,4)$ case, we do not enumerate the candidates for 
$c_1$ and $c_1^{(2)}$, but use the upper bounds for $c_1$ from (47) and (49) 
and the lower bounds for $c_1^{(2)}$ from (51) and (53) in (56). Thus using 
(51), (52) and (56), the Basic Method yields candidates for $c_1^{(1)}$ which 
depend on each linked pair $(c_0^{(1)}, c_0^{(2)})$ (We drop $c_0$). Note that, 
from (46) 
$$ \beta_1^2 \mid \beta_1 c_1^{(1)} + c_0^{(1)} - c_0^{(2)} , $$
and we further reduce our list of candidates to satisfy this condition. The 
total number of triples $(c_1^{(1)}, c_0^{(1)}, c_0^{(2)})$ at this stage 
is 2071.

Again  using the Basic Methods, we determine, independently of the 
foregoing calculations, the candidates for $c_3$ and $c_3^{(1)}$ (there are 452
and 187 respectively). Now 
\begin{equation}
c_2^{(1)} = \frac{1}{2}(c_3 + \beta_1 c_3^{(1)} + {c_3^{(1)}}^2).
\end{equation}
The basic inequalities that $c_2^{(1)}$ and its conjugates must satisfy 
together with the fact  that the second derivative of 
$p_1^{\sigma_i}(x), i = 2,3$ must have two real roots in the interval
$(\beta_i, 0)$ gives inequalities relating $c_3^{(1)} $ and $c_2^{(1)}$ and 
hence involving $c_3$ and $c_3^{(1)}$. We thus determine all pairs 
$(c_3, c_3^{(1)})$ which are linked by these inequalities. Furthermore, 
since $2 \mid c_3 + \beta_1 c_3^{(1)} + {c_3^{(1)}}^2$ we reduce the pairs 
to satisfy this divisibility condition. We then solve for $c_2^{(1)}$ using (57) and drop $c_3$. There are 2218 resulting pairs $(c_3^{(1)}, c_2^{(1)})$.

We now relate these linked pairs $(c_3^{(1)}, c_2^{(1)})$ to the linked pairs 
$(c_0^{(1)}, c_0^{(2)})$ by inequalities. Once again using the AM/GM 
inequality yields $\sigma_i(c_3^{(1)}) \geq 4 \sigma_i(c_0^{(1)})^{1/4}$
and $\sigma_i(c_2^{(1)}) \geq 6 \sigma_i(c_0^{(1)})^{1/2}$ for $i=2,3$.
Also $\sigma_i(c_2^{(1)}) < 3 \sigma_i(c_0^{(1)})^{1/2} +  3\beta_i^2$
for $i = 2,3$. A bit of manipulation using the AM/GM inequality also yields 
$c_2^{(1)} > 2 \sqrt{3} {c_0^{(1)}}^{1/2}$. We thus determine all 4-tuples which
are linked by these inequalities. There are a total of 74570.

We now have a collection of 1051111 5-tuples $(c_3^{(1)}, c_2^{(1)}, c_1^{(1)}, c_0^{(1)}, c_0^{(2)})$ indexed by the linked pairs $(c_0^{(1)}, c_0^{(2)})$.
They must satisfy equation (46). Implementing this gives 1934 4-tuples (we drop
$c_0^{(2)})$. Then requiring that the first derivative of $p_1^{\sigma_i}(x)$,
$i = 2,3$ has  three roots in the interval $(\beta_i,0)$ (see \S 7) yields 
a list of 746 polynomials. These and their conjugates can then be 
numerically solved and only 8 polynomials have the correct distribution of 
real roots. All these 8 turn out to be reducible and so we do not obtain 
any groups in this $(7,7,4)$ case.

\noindent {\bf Comments on the other five cases}

Case $(8,8,5)$. This is tackled in a very similar manner to the preceding
$(7,7,4)$ case. The main difference is that in this case, we use the 
inequalities (47) to (56) to enumerate the candidates for $(c_1, c_1^{(1)}, c_1^{(2)})$
depending on the linked triple $(c_0, c_0^{(1)}, c_0^{(2)})$ which satisfy 
equation (55).  As in the preceding case, we then determine the pairs 
$(c_4^{(1)}, c_3^{(1)})$ and relate them by inequalities to the linked pair $(c_0^{(1)}, c_0^{(2)})$.
In this case
\begin{equation}
\beta_1^5 + c_4^{(1)} \beta_1^4 + c_3^{(1)} \beta_1^3 + c_2^{(1)} \beta_1^2 + c_1^{(1)} \beta_1 + c_0^{(1)} = - c_0^{(2)}
\end{equation}
\begin{equation}
5 \beta_1^4 + 4 c_4^{(1)} \beta_1^3 + 3 c_3^{(1)} \beta_1^2 + 2 c_2^{(1)} \beta_1 + c_1^{(1)} = c_1^{(2)}
\end{equation}
from which we obtain
\begin{equation}
3 \beta_1^5 + 2 c_4^{(1)} \beta_1^4 + c_3^{(1)} \beta_1^3 - c_1^{(1)} \beta_1 - c_1^{(2)} \beta_1 - 2 c_0^{(1)} - 2 c_0^{(2)} = 0.
\end{equation}
We now determine all  6-tuples $(c_4^{(1)}, c_3^{(1)}, c_1^{(1)}, c_1^{(2)}, c_0^{(1)}, c_0^{(2)})$ which satisfy (60) and use (58) to determine 
$c_2^{(1)}$ from the remaining coefficients. Now as for degree 3, we reduce
our collection by the condition that, at the non-identity place, the degree 
3 polynomial which is the second derivative of $p_1(x)$ has three real roots 
in the interval $(\beta_2, 0)$ where $\beta_2 = -(2 + \sqrt{2})$. This gives 
us 95 polynomials which can then be solved numerically and none have five 
real roots in the interval $(\beta_2, 0)$. So there are no groups in this case.

Case $(8,8,4)$. A simplified version of the above yields three polynomials 
with the correct numbers of real roots, but at the identity embedding, the real roots do not lie in the interval $(\beta_1, 0)$ where $\beta_1 = -(2-\sqrt{2})$.

Case $(12,12,4)$. In this case, five polynomials have the correct distribution of roots, but for 3 of them, the real roots do not lie in the interval 
$(\beta_1, 0)$ at the identitiy embedding and for the other two, the resulting
$\gamma(12,12)$ value lies outside the contour $\Omega_{12,12}$.

Case $(6,6,5)$. An even more simplified version of the above method yields 31 
polynomials with 3 real roots in the interval $(-1,0)$. For all but one of them,
the associated $\gamma(6,6)$ lies outside the contour $\Omega_{6,6}$. Thus there
is one candidate to be considered by geometric methods.

Case $(6,6,4)$. Using the same techniques as above, there are 70 polynomials
with the correct distribution of roots and such that $\gamma(6,6)$ lies inside
the contour $\Omega_{6,6}$, all of these needing further examination by 
geometric methods.

  \section{Finite Co-volume}

  In \S\S 5, 6, 7 and 8, we have outlined the methods we applied to all the triples
  $(p,q,r)$ which appear on the {\it Aspiring List}. The result is a set
  of irreducible polynomials of degree $r$ over $L$ whose complex
  roots $\gamma$ satisfy all four conditions (alternatively Lemma \ref{condition4'}) of 
  Theorem \ref{2genthm} and also the inequalities of \S 3,  meaning they are not obviously of infinite volume. 
  
  This means that $\gamma$
  determines a group $\G$ generated by elements of orders $p,q$ which is
  a subgroup of an arithmetic Kleinian group and hence discrete.   From a specific value of $\gamma$ we can compute 
  the (normalized) matrices $A$ and $B$ which represent the generators
  $f,g$ of $\G$ (see \S 3). The inequalities of \S 3 are derived from the 
  geometric result that, if $\G$ is to be of finite co-volume then it cannot be 
  a free product so that the isometric circles of $g$ and $g^{-1}$ 
  cannot lie within the intersection of the isometric cicles of $f$ and
  $f^{-1}$. 
  
  Two further, but more complicated, geometric conditions
  (temporarily labelled {\it Free2} and {\it Free3} for use in the Examples 
  below), necessary for $\G$ not to be a free product have been given in 
  \cite{MM2} in terms of the locations of the isometric circles
  of combinations of $f,g$. These simply consist of looking at the images of 
  the isometric circles of one generator,  say $g$ under the transformation 
  $f$ and trying to piece together a fundamental domain from the intersection 
  pattern.  For instance, illustrated below,  although the isometric circles 
  of $g$ do not lie in the region bounded between the isometric circles of 
  $f$ we have $f(I(g)\cup I(g^{-1}))\cap ( I(g)\cup I(g^{-1})) =\emptyset$ 
  and so if we look at the region (where we write $I(g)$ to mean the 
  disk bounded by $I(g)$ etc)
  \[ I(f)\cap I(f^{-1}) \setminus (I(g)\cup I(g^{-1}) \cup f(I(g)\cup I(g^{-1})) \]
shaded below one can show without too much effort that this region lies within a fundamental domain for $\langle f,g \rangle$ on $\hat{\IC}$ and in fact $\langle f,g \rangle$ is free on generators.

\bigskip
  
 \hskip30pt\scalebox{0.60}[0.60]{ \includegraphics*[viewport=80 200 500 440]{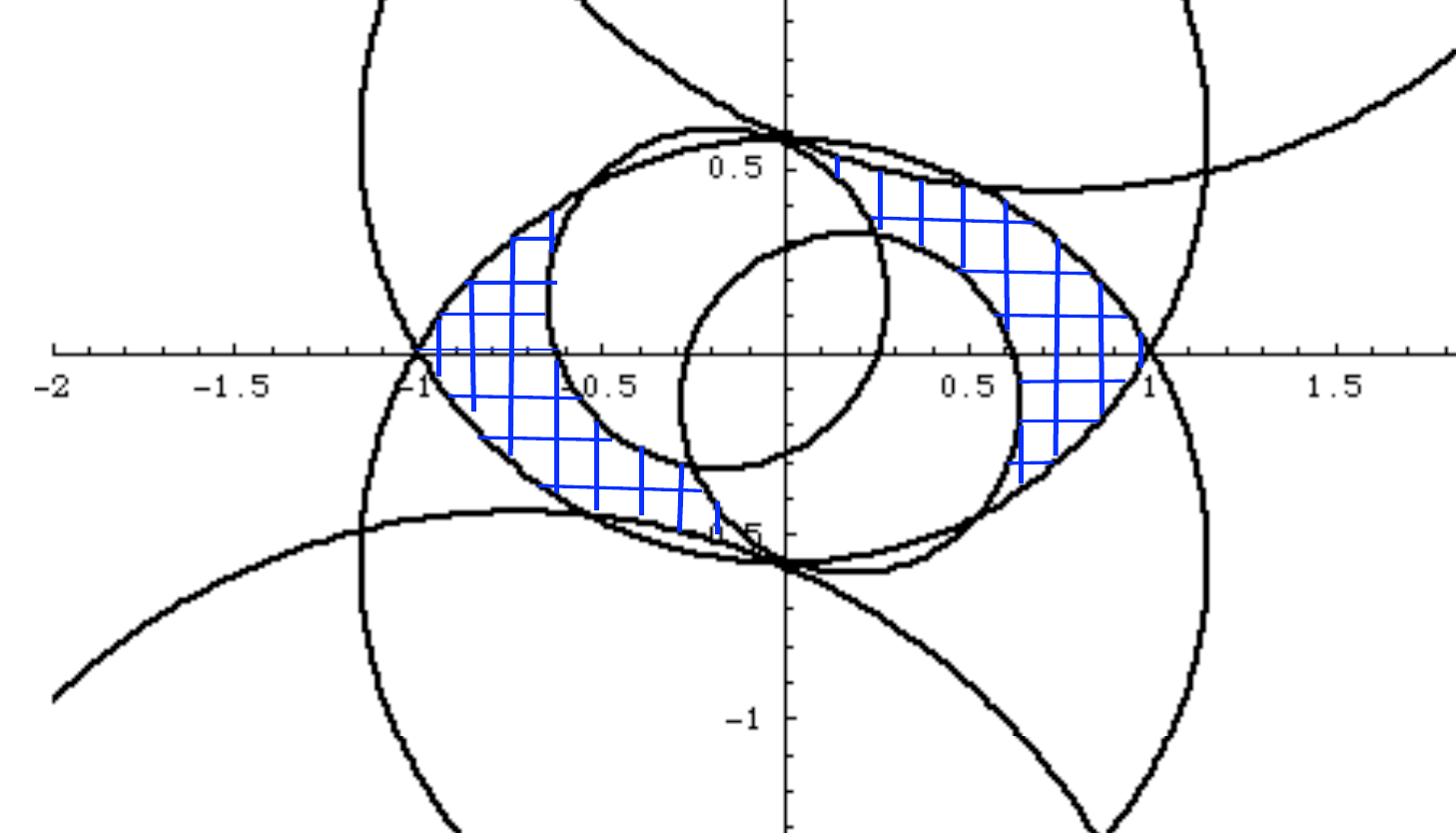}}~
\begin{center}
 Diagram 3. Second level isometric circles
  \end{center}
  
 \bigskip
 Of course one can go on looking at more and more isometric circles and their 
patterns and using this to formulate algebraic inequalities on $\gamma$.  
However after three levels this becomes quite impractical and we shall 
discuss below the computer program used to deal with these cases.
 
  We apply these two elementary tests on $\gamma$ to reduce the list of polynomials
arising from \S 6,7 and 8.

\medskip

\begin{example}
$(p,q,r) = (10,10,3)$. From Table 8, there are 5 candidates for $\gamma$ which satisfy 
the four conditions of  Theorem \ref{2genthm} and the inequalities of \S 3. Their minimum
polynomials over $L = \IQ(\sqrt{5})$ are given below.
\begin{center}
\begin{tabular}{|c|c|c|} \hline
No. & $\gamma$ & polynomial \\ \hline
 1 & $-4.918226 + 5.698268 i$ & $x^3 +\frac{13+3\sqrt{5}}{2}x^2+(30+12\sqrt{5})x+1 $ \\ \hline
2 & $0.635991 + 5.238279 i$ & $x^3 + (1-\sqrt{5})x^2+\frac{31+11\sqrt{5}}{2}x+1$ \\ \hline
3 & $3.251943 + 8.478242 i $ & $x^3 + (-2-2\sqrt{5})x^2+(42+18\sqrt{5})x+\frac{3+\sqrt{5}}{2} $ \\ \hline
4 & $ 8.794158 + 4.828433 i $ & $x^3+\frac{-15 -9\sqrt{5}}{2}x^2+(51+22\sqrt{5})x+\frac{3+\sqrt{5}}{2} $ \\ \hline
5 & $ 6.180432 + 10.631111 i $ & $x^3+\frac{-9-7\sqrt{5}}{2}x^2+(77+33\sqrt{5})x+(3+\sqrt{5}) $ \\ \hline
\end{tabular}
\end{center}
Test {\it Free2} removes cases 1 and 5, while test ${\it Free3}$
removes case 2 and 3, leaving just one possibility to consider further.
\end{example}

 \hskip30pt\scalebox{0.50}[0.50]{ \includegraphics*[viewport=80 120 600 520]{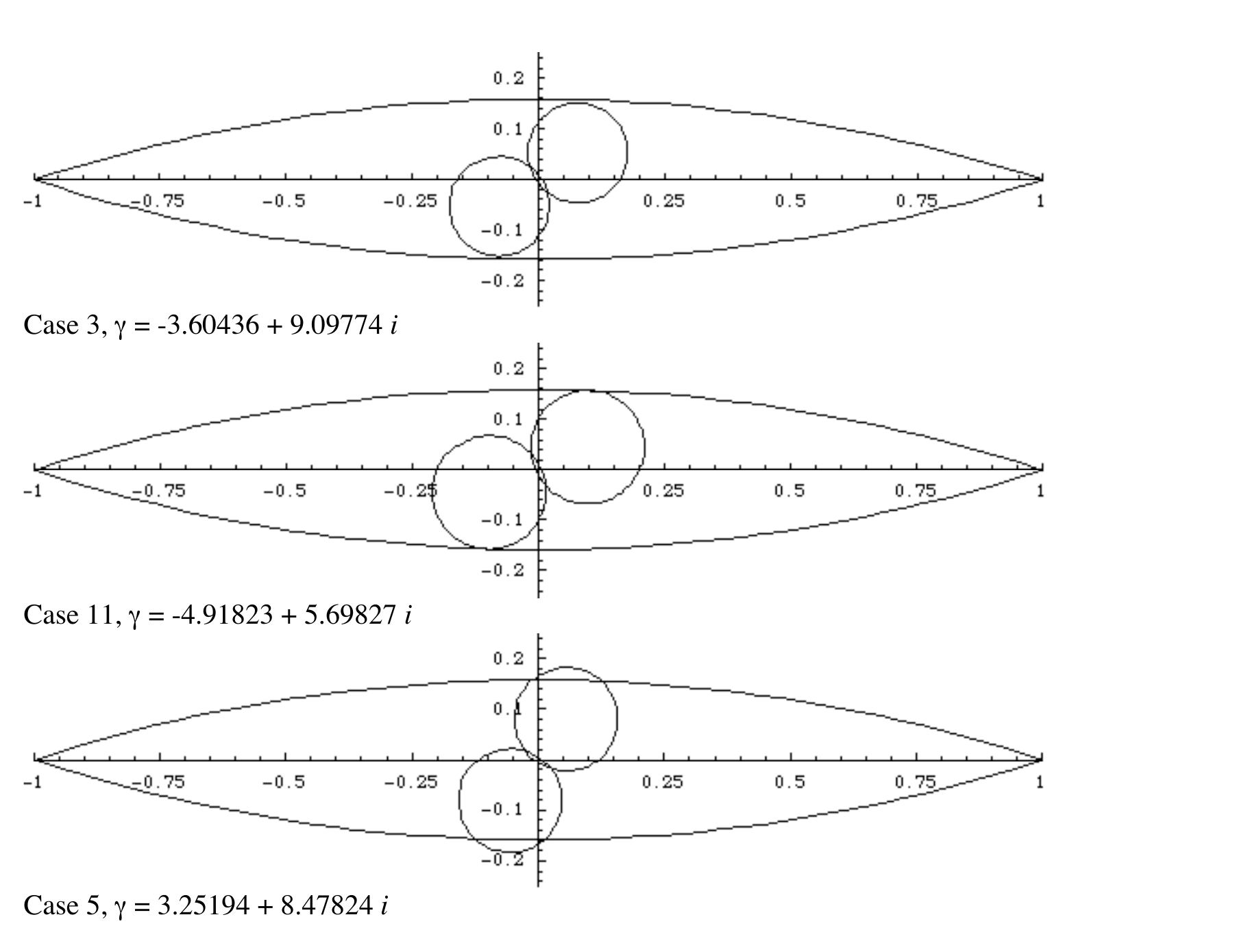}}
\begin{center}
Diagram 4. $(10,10,3)$-cases.
\end{center}
\medskip
We list the numbers of  candidates brought forward from \S 6,7 and 8 and in 
Table 9 show that many are eliminated using these Tests. Thus under ``FT2'',
``FT3'',
are the number eliminated using FreeTest2 and FreeTest3 respectively. We also 
remove any reducible polynomials or duplicates which have survived to this 
stage. The final
column gives the numbers of polynomials which are passed to the next step in 
procedure given below.
\begin{table}[h]
\begin{center}
\begin{tabular}{|c|c|c|c|c|} \hline
Triple & No. & FT2 & FT3 & Rem. \\ \hline
$(30,30,2)$ & 1 & 0 & 0 & 1 \\ \hline
$(18,18,2)$ & 2 & 1 & 1 & 1 \\ \hline
$(18,9,2)$ & 2 & 2 & 0 & 0 \\ \hline
$(14,7,2)$ & 1 & 0 & 1 & 0 \\ \hline
$(12,12,2)$ & 18 & 7 & 4 & 7 \\ \hline
$(12,6,2)$ & 1 & 0 & 1 & 0 \\ \hline
$(10,10,2)$ & 15 & 8 & 2 & 5 \\ \hline
$(9,9,2)$ & 2 & 1 & 1 & 0 \\ \hline
$(8,8,2)$ & 20 & 8 & 4 & 8 \\ \hline
$(7,7,2)$ & 8 & 4 & 2 & 2 \\ \hline
$(6,6,2)$ & 24 & 9 & 5 & 10 \\ \hline
$(12,12,3)$ & 19 & 10 & 5 & 4 \\ \hline
$(12,6,3)$ & 1 & 1 & 0 & 0 \\ \hline
$(10,10,3)$ & 5 & 2 & 2 & 1 \\ \hline
$(8,8,3)$ & 29 & 13 & 10 & 5 \\ \hline
$(7,7,3)$ & 1 & 0 & 1 & 0 \\ \hline
$(6,6,3)$ & 124 & 52 & 30 & 38 \\ \hline
$(6,6,5)$ & 1 & 0 & 0 & 1 \\ \hline
$(6,6,4)$ & 70 & 36 & 19 &15 \\ \hline
\end{tabular}
\caption{Geometric Test Results}
\end{center}
\end{table}

Finally, a computer program has been developed, initially by J. McKenzie, and named JSnap to study subgroups $\G$ of $\PSL(2,\IC)$ which have two generators of finite order. This is effectively an implementation of the Dirichlet routine in J. Weeks' program Snappea.  This program aims to find a Dirichlet region for the group $\G=\langle f, g\rangle$.  A very important point to note here is that we know {\em a priori} that the group in question is discrete.  However it is theoretically possible that the group $\G$  is geometrically infinite and so computationally impossible to identify a fundamental domain.  JSnap runs and either produces a fundamental domain - either of finite or infinite volume - or produces an error message if it can't put together a fundamental domain after looking at words of a given bounded length.  In our situation JSnap always produces a fundamental domain which is either compact or meets the sphere at infinity in an open set (which itself will be a fundamental domain for the action of $\G$ on $\hat{\IC}\setminus \Lambda(\G)$).  In this latter case the group cannot be of finite co-volume (and it might also not be free on generators - for instance certain Web-groups may arise) and so we can eliminate these cases.

If the fundamental domain found by JSnap is compact,  then JSnap also returns an approximate co-volume.

\medskip

  In this way the  remaining possibility for   $(p,q,r)=(10,10,3)$  
$\gamma = 8.794158 + 4.828433 i$
 is shown to have a fundamental domain which meets the sphere at infinity in 
an open set,  and thus cannot be arithmetic.  

 Applying  JSnap to the 15 cases in $(6,6,4)$ shows
that they all have infinite volume and so there are no corresponding
arithmetic Kleinian groups.  It is a similar story for $r=3$ except that here we meet our first example whose isometric circle configuration and $\gamma$ value are illustrated below.

 \hskip30pt\scalebox{0.85}[0.85]{ \includegraphics*[viewport=90 110 500 350]{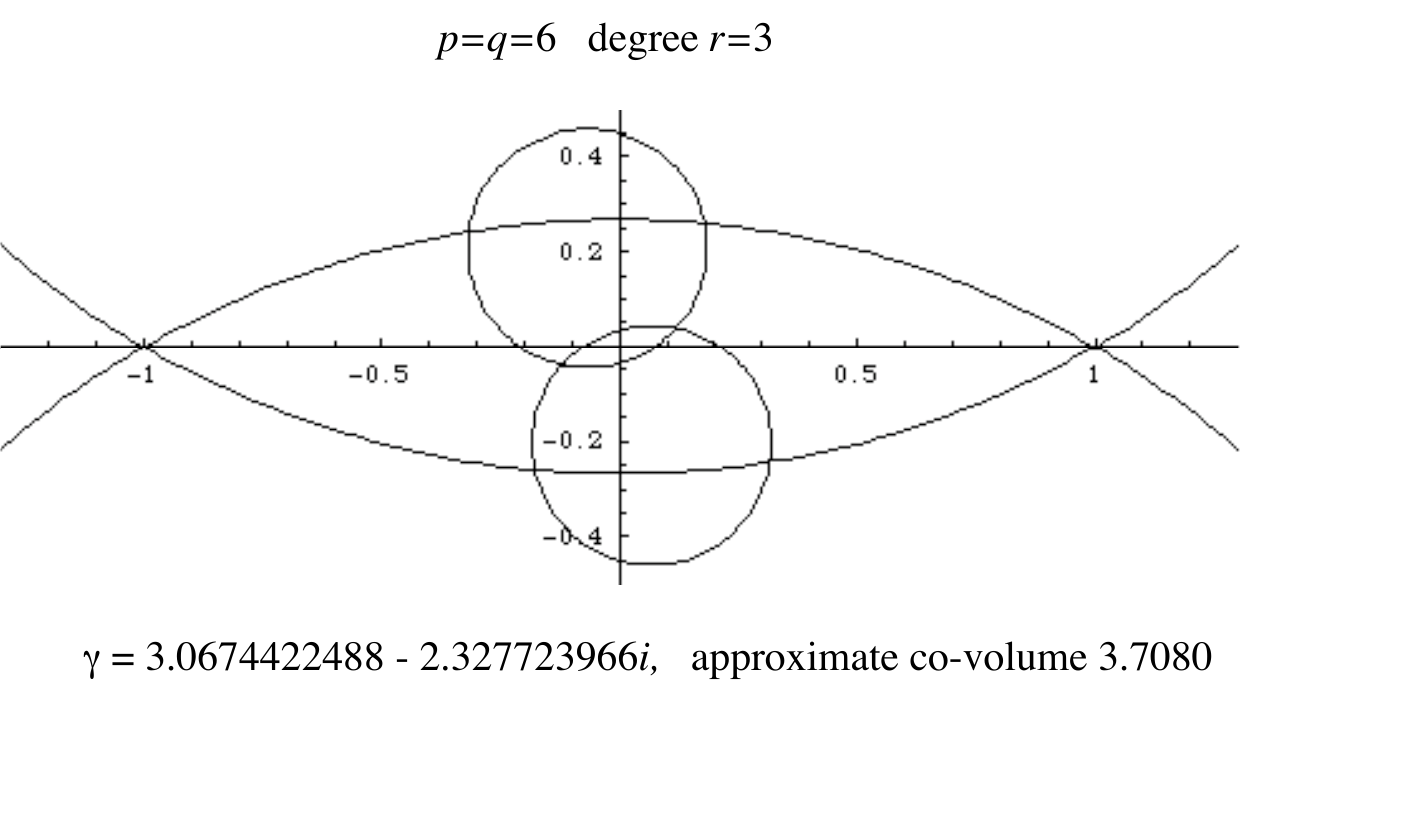}}
\begin{center}
Diagram 5. A finite co-volume $(6,6,3)$ example.
\end{center}

\section{The end product}

After dealing with all the cases on the Aspiring List in the manner outlined in the proceeding sections we are left with just 17 complex values of $\gamma$ whose corresponding group 
 JSnap identifies as having finite co-volume and giving us an approximation to this co-volume. The groups corresponding to 
all but one of these
values have already been identified as arithmetic in the 
literature in \cite{HLM1} as being obtained from surgery on 2-bridge
knots and links. This can also be ascertained using Snappea as discussed in the
introduction. The arithmetic data required to define these groups can be
recovered from the polynomials satisfied by $\gamma$. 
The value $\sqrt{-3}$ gives a group which is not co-compact 
and has been discussed in 
\cite{MM2}. In addition there is one further group, also not co-compact,
which corresponds to the only real value of $\gamma$ which arises for 
groups with generators of orders $\geq 6$ (see \cite{MM3,MM4}).
 
 \bigskip

In the tables below we list all the data on the groups we have found.
The polynomial is that satisfied by $\gamma=\gamma(f,g)$ where $o(f)=p,
o(g)=q$ over the field $\IQ(\cos 2 \pi/p, \cos 2 \pi/q)$. In most cases the
description is given as an orbifold obtained by surgery on the one 
boundary component of a 2-bridge knot or on both boundary components
of a 2-bridge link.

 The appearance of the same description occurring 
twice in this table is discussed in the introduction as identification of different Nielsen classes of generators. In the other two cases, the description
refers to \cite{MM2}. 

The commensurability class of the arithmetic group
is, as we have discussed, determined by the field of definition and the
defining quaternion algebra. In all the cases here, the discriminant
$\Delta$ in the table uniquely describes the field given its degree and 
that it has one complex place. The quaternion algebra is determined by its 
ramification set which must include all real places so that the finite
ramification suffices to identify the quaternion algebra. (The convention
used here is that if a rational prime $p$ splits as ${\cal P}_p {\cal P}_p'$
then these are ordered so that $N({\cal P}_p) \leq N({\cal P}_p')$.)
Note that the only commensurable pairs are those co-compact 
groups which are actually
equal and given by different Nielsen classes of generators and the non-co-compact
pair.

The orbifold volume is determined by Snappea and JSnap and the minimum
volume is the smallest volume of any orbifold in the commensurability
class which can be determined from the arithmetic data defining
the commensurability class \cite{Bo,MR}.
 
\begin{center}
\begin{tabular}{|c|c|c|c|} \hline
No. & $(p,q)$ & $\gamma$ & poly  \\ \hline
1 & $(12,12)$ &$-0.259113+1.998874 i$ & $x^3+(4-2\sqrt{3})x^2+(11-4\sqrt{3})x+(7-4\sqrt{3})$  \\ \hline
2&$(12,12)$ & $-0.633975+0.930605i$ & $x^2+(3-\sqrt{3})x+(3-\sqrt{3})$ \\ \hline
3 & $(10,10)$ & $-1 +2.058171i$ & $x^2+2x + (3+\sqrt{5})$ \\ \hline
4 & $(8,8)$ & $-0.792893+0.978318i$ & $x^2+(3-\sqrt{2})x + (3-\sqrt{2}) $ \\ \hline
5 & $(6,6)$ & $-1.877438+0.744861i$ & $x^3+4x^2+5x+1$ \\ \hline
6 & $(6,6)$ & $-2.884646+0.589742i$ & $x^3+6x^2+10x+2$ \\ \hline
7 & $(6,6)$ & $-0.891622+1.954093i$ & $x^3+2x^2+5x+1 $ \\ \hline
8 & $(6,6)$ & $1.092519+2.052003i$ & $x^3-2x^2+5x+1 $ \\ \hline
9 & $(6,6)$ & $3.067442+2.327724i$ & $x^3-6x^2+14x+2 $ \\ \hline
10 & $(6,6)$ & $0.124046+2.836576i$ & $x^3+8x+2$ \\ \hline
11 & $(6,6)$ & $2.124407+2.746645i$ & $x^3-4x^2+11x+3$ \\ \hline
12 & $(6,6)$ & $4.109638+2.431700i$ & $x^3-8x^2+21x+5$ \\ \hline
13 & $(6,6)$ & $-1+i$ & $x^2+2x+2$ \\ \hline
14 & $(6,6)$ & $-2+1.414214i$ & $x^2+4x+6$ \\ \hline
15 & $(6,6)$ & $1.732051i$ & $ x^2+3 $ \\ \hline
16 & $(6,6)$ & $-1+2.645751i$ & $x^2+2x+8$ \\ \hline
17 & $(6,6)$ & $1+3i$ & $x^2-2x+10$ \\ \hline
18 & $(6,6)$ & $-1$ & $x+1$ \\ \hline
\end{tabular}
\\[\baselineskip]
\begin{tabular}{|c|c|c|c|c|c|} \hline
No. & Description & $\Delta$ & ${\rm Ram}_f(A)$ & Orb. Vol. & Min. Vol. \\ \hline
1 & $(12,0),(12,0)$ on $8/3$ & $-288576$ & $\emptyset$ & $3.3933$ & $0.424167$ \\ \hline
2 & $(12,0)$ on $5/3$ & $-1728$ & ${\cal P}_2, {\cal P}_3$ & $ 1.8026$ & $0.450658$ \\ \hline
3 & $(10,0)$ on $13/5$ & $-400$ & ${\cal P}_2, {\cal P}_5$ & $5.1674$ & $1.291862$ \\ \hline
4 & $(8,0)$ on $5/3$ & $-448$ & ${\cal P}_2, {\cal P}_7$ & $1.5438$ & $0.385966$ \\ \hline
5 & $(6,0)$ on $7/3$ & $-23$ & ${\cal P}_3$ & $2.0425$ & $0.510633$ \\ \hline
6 & $(6,0),(6,0)$ on $20/9$ & $-76$ & ${\cal P}_2,{\cal P}_3,{\cal P}_3'$ & $5.2937$ & $0.661715$ \\ \hline
7 & $(6,0),(6,0)$ on $8/3$ & $-31$ & ${\cal P}_3$ & $2.6386$ & $0.065965$ \\ \hline
8 & $(6,0)$ on $7/3$ & $-23$ & ${\cal P}_3$ & $2.0425$  & $0.510633$ \\ \hline
9 & $(6,0)$ on $13/3$ & $-44$  & ${\cal P}_2$ & $3.7068$ & $0.066194$ \\ \hline
10 & $(6,0)$ on $13/3$ & $-44$ & ${\cal P}_2$ & $3.7068$ & $0.066194$ \\ \hline
11 & $(6,0)$ on $15/11$ & $-31$ & ${\cal P}_3'$ & $4.2217$ & $0.263861$ \\ \hline
12 & $(6,0)$ on $65/51$ & $-23$ & ${\cal P}_5$ & $8.7986$ & $0.078559$ \\ \hline
13 & $(6,0)$ on $5/3$ & $-4$ & ${\cal P}_2, {\cal P}_3$ & $1.2212$ & $0.305322$ \\ \hline
14 & $(6,0),(6,0)$ on $12/5$ & $-8$ & ${\cal P}_2, {\cal P}_3$ & $4.0153$ & $0.250960$ \\ \hline
15 & Non-compact  $\Gamma_{21}$  & $-3$ & $\emptyset$ & $1.0149$ & $0.253735$ \\ \hline
16 & $(6,0),(6,0)$ on $30/11$ & $-7$ & ${\cal P}_2, {\cal P}_3$ & $7.1113$ & $0.888915$ \\ \hline
17 & $(6,0),(6,0)$ on $24/7$ & $-4$ & ${\cal P}_2, {\cal P}_5$ & $6.1064$ & $0.152661$ \\ \hline
18 & Non-compact $\Gamma_{20}$ & $-3$ & $\emptyset$ & $0.5074$ & $0.084578$ \\ \hline 
\end{tabular}
\end{center}

\bigskip 

\noindent {\bf Remark.}  In eliminating certain candidates because they are not of co-finite volume - and which are guaranteed discrete by our arithmetic criteria - we ran into a number of interesting examples where our computational package JSnap had difficulty.   This was largely to do with accumulation of roundoff error.  In exploring these groups (looking for regions of discontinuity) we made use of the package``lim" developed by C. McMullen to draw the limit sets of Kleinian groups.  The limit set of one such group is illustrated below.  After seeing these pictures we were encouraged to modify our code to run on a different platform with higher precision to get an infinite volume fundamental region.  However it is clear that in these sorts of cases (with parameters algebraic integers of low degree) that working with a version of Snap (the precise arithmetic version of Snappea developed by O. Goodman et al, \cite{asnap}) would be the correct way forward.  We are currently developing this program which will surely be necessary in extending our results beyond the cases $p,q\geq 6$.

 \scalebox{0.4}{\includegraphics*[viewport=-150 100 1000 1000]{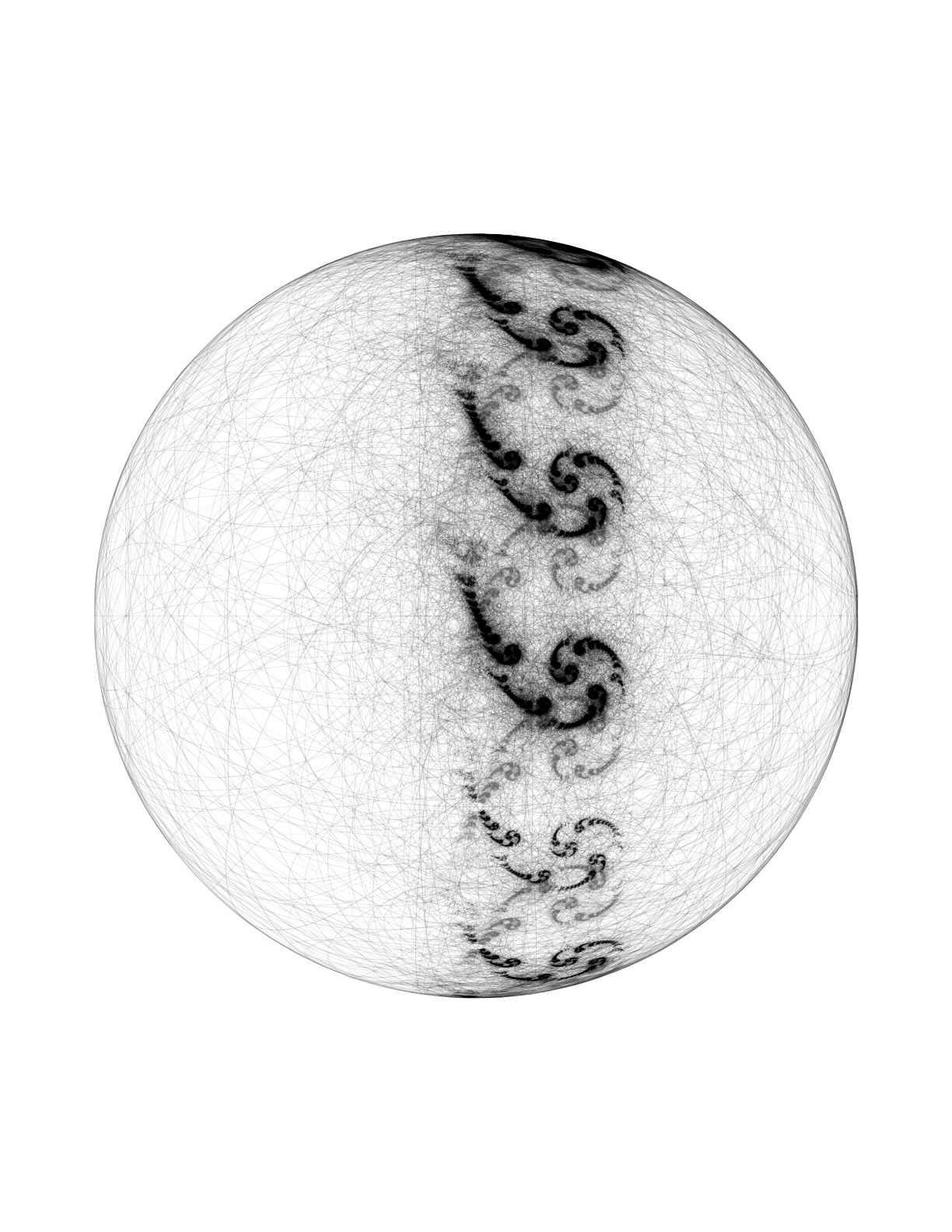}}
\begin{center}
Limit set of Kleinian group with two generators of order 12 and  $\gamma= 2.73205+3.193141i$,  a root of $x^2-2(1+\sqrt{3})x+(9+5\sqrt{3})=0$ 
\end{center}

\subsection{Generalised Triangle Groups}\label{proofngtg}

Here we prove Corollary \ref{notrigrp}.  First we note that the only groups we need to consider here are the surgeries on two-bridge knot and link groups and thus the following lemma will suffice.

\begin{lemma} Let $\G$ be the orbifold fundamental group of $(p,0)$-$(q,0)$ ($p,q\geq 2$) Dehn surgery on a two-bridge knot or link.  Then $\G$ does not have a presentation as a generalised triangle group.
\end{lemma}
\noindent{\bf Proof.}  Every element $g$ of finite order in $\G$ has a nontrivial fixed point set in $\IH^3$ which projects to an edge in the singular set of $\IH^3/\G$.  Elements in the same conjugacy class project to the same edge.  Next $(p,0)$-$(q,0)$ Dehn surgery on a two-bridge knot or link has at most two components in its singular set (one if it is a knot).  Let us denote the two primitive generators of order $p$ and $q$ arising from this surgery as $f$ and $g$,  so $\G=\langle f,g\rangle$.   Suppose $\G$ has a presentation of the form 
\begin{equation}\label{gtp}
\langle a^r=b^s=w(ab)^t\rangle, \hskip20pt r,s,t \geq 2
\end{equation}
  There are at most two conjugacy classes of torsion in $\G$.  Thus $a, b$ and any other element of finite order are conjugates of elements of $\langle f \rangle $ or $\langle g \rangle$.    
  Thus,  possibly increasing $r$ or $s$ and the complexity of $w$,  we see that we can find a presentation of the form (\ref{gtp}) with $a$ and $b$  conjugates of $f$ and $g$ and so  $\{r,s\}\subset \{p,q\}$.  
  
  Suppose that $b$ is not a conjugate of $a$.  Then as $w$ must also be conjugate into $\langle f \rangle $ or $\langle g \rangle$, the relation $w^r=1$ is a direct consequence of the relators $a^p=b^q=1$ and so $\G= \langle a \rangle \ast \langle b \rangle$ which is not possible for a co-finite volume lattice.  Thus $b$ is a conjugate of $a$ and $r=s$.  This quickly implies that the abelianisation of $\G$ is a subgroup of $\langle a \rangle$ as $w$ reduces to a power of $a$.  Thus, if $\G$ has a presentation as at (\ref{gtp}) we have deduced that $\G$ abelianises to a cyclic group $\IZ_k$ with $k|p$ or $k|q$. Further,  we cannot be dealing with a knot surgery  as $a$ not conjugate to $w$ implies two components to the singular set.  
  
  Next, every two bridge link has a presentation on a pair of meridians of the form $\langle u,v: u w= w u\rangle$ for $w$ a word determined by the Schubert normal form \cite{BZ}. Dehn surgery is equivalent to adding the relators $u^p=v^q=1$ which then gives $\IZ_p+\IZ_p$ as the abelianisation.  
  
  Thus there can be no presentation as at (\ref{gtp}) and the proof of the lemma is complete.\hfill $\Box$

\end{document}